\documentclass{amsart}
\usepackage{amssymb,amsthm,amsmath,amsxtra}

\newtheorem{theorem}{Theorem}[section]
\newtheorem{df}[theorem]{Definition}
\newtheorem{lem}[theorem]{Lemma}
\newtheorem{cor}[theorem]{Corollary}
\newtheorem{prop}[theorem]{Proposition}

\newtheorem*{namedtheorem}{\theoremname}
\newcommand{\theoremname}{testing}

\theoremstyle{definition}

\newtheorem{remark}[theorem]{Remark}

\newcommand{\h}{\mathfrak H}
\renewcommand{\O}{\mathcal O}
\renewcommand{\a}{\mathfrak a}
\newcommand{\Af}{\mathbb A_f}
\newcommand{\qq}{\mathbb Q}
\newcommand{\zz}{\mathbb Z}

\title{Borcherds forms and generalizations of singular moduli}
\author{Jarad Schofer}
\begin{document}

\maketitle

\section{Introduction}

Borcherds forms are meromorphic modular forms for arithmetic subgroups
$\Gamma$ of the orthogonal group
$O(n,2)$ which arise from a regularized theta lift of (vector valued)
modular forms of weight $1-\frac{n}{2}$
for $SL_2(\zz)$ with poles at the cusp.  They have interesting product expansions and
explicitly known divisors (cf. \cite{Bor98}).  In some cases they can be
realized as classical modular forms, such as the difference of two modular
$j$-functions or as the discriminant function $\Delta$ (see \cite{Bru}).  In this paper, we give a factorization of values of
Borcherds forms at CM points.  The main
result can be viewed as a generalization of the singular moduli
result (Theorem 1.3 of \cite{GZ}) of Gross and Zagier. In fact, our method gives a
new proof of their result, which will be discussed in
a sequel to this paper.

\par Let $V$ be a vector space with quadratic form $Q$ of
signature $(n,2)$ and let $D$ be the space of oriented negative-definite
two-planes in $V(\mathbb R)$.  $D$ is the symmetric space for $O(n,2)$
and has a Hermitian structure.  For example, when $n=1$,
$D\simeq\h^+\sqcup\h^-$ is the union of the upper and lower
half-planes $\h=\h^+$ and $\h^-$, respectively.  Let
$H=\text{GSpin}(V)$ and let $K\subset H(\Af)$ be a compact open
subgroup, where $\Af$ is the finite adeles.  We consider the
quasi-projective variety
\[
X_K=H(\qq)\backslash\Big(D\times H(\Af)/K\Big)\simeq\coprod_j\Gamma_j\backslash D^+,
\]
for a finite number of arithmetic subgroups $\Gamma_j\subset H(\qq)$, and where $D^+\subset D$ is
the subset of positively oriented two-planes. 

\par Recall the theory of Borcherds forms on $X_K$.  For a lattice
$L\subset V$ with dual
\[
L^\vee=\{x\in V\mid (x,L)\subseteq\zz\}
\]
such that $L^\vee\supset L$, there exists a finite dimensional
subspace $S_L\subset S(V(\Af))$ of the Schwartz space of $V(\Af)$
defined as follows.  Let $\hat{L}=L\otimes_\zz\hat{\zz}$.  Then $S_L$
is the space of functions with support in $\hat{L}^\vee$ which are constant
on cosets of $\hat{L}$.  A natural basis of $S_L$ is
\[
\{\varphi_\eta=\text{char}(\eta+L)\mid\eta\in L^\vee/L\}
\]
and dim $S_L=| L^\vee/L|$.  There exists a representation
$\omega$ of (the metaplectic extension $\Gamma'$ of)
$\Gamma=SL_2(\zz)$ on $S(V(\Af))$ preserving $S_L$; see section 4 of
\cite{Bor98} for details.

\par A modular form $F:\h\to S_L$ of weight
$1-\frac{n}{2}$ and type $\omega$ for $\Gamma$ satisfies
\[
F(\gamma\tau)=(c\tau+d)^{1-\frac{n}{2}}\omega(\gamma)(F(\tau))
\]
for all $\gamma=\left(\begin{smallmatrix}
a&b\\
c&d
\end{smallmatrix}\right)\in\Gamma$.  The function $F$ has Fourier expansion
\begin{equation}\label{introFourier}
F(\tau)=\sum_\eta F_\eta(\tau)\varphi_\eta=\sum_\eta\sum_m
c_\eta(m)\mathbf q^m\varphi_\eta,
\end{equation}
where $\mathbf q=e^{2\pi i\tau}$.  We say that $F$ is \emph{weakly holomorphic} if only a finite number of the $c_\eta(m)$'s with $m<0$
are non-zero.  Furthermore, we call such a modular form
\emph{integral} if the non-positive Fourier coefficients lie in $\zz$.

\par To a weakly holomorphic integral modular form $F$ of weight $1-\frac{n}{2}$, Borcherds attaches a
function $\Psi(F)$ (called a Borcherds form), which is a meromorphic
modular form on the space $D\times H(\Af)$ with respect to $H(\qq)$.  The weight of $\Psi(F)$ is
$\frac{1}{2}c_0(0)$ and the divisor of $\Psi(F)^2$ is given explicitly in terms of
the negative Fourier coefficients of $F$,
\[
\text{div}(\Psi(F)^2)=\sum_\eta\sum_{m>0}c_\eta(-m)Z(m,\eta,K),
\]
for divisors $Z(m,\eta,K)$ on $X_K$.  The concrete connection
between $F$ and $\Psi(F)$ is given by a regularized theta lift
\[
\Phi(z,h;F):=\int_{\Gamma\backslash\h}^{\bullet}((\;F(\tau),\theta(\tau,z,h)\;))v^{-2}dudv,
\]
where $z\in D, h\in H(\Af)$ and $\tau=u+iv\in\h$, and where
$((\;F(\tau),\theta(\tau,z,h)\;))$ is a theta function constructed
from the Fourier expansion of $F$; see section \ref{basicsetup} for
details.  Since $F$ has a pole at the cusp, this integral diverges and
so it must be regularized.  See \cite{Bor98} or section
\ref{basicsetup} for the exact definition of the regularized
integral.  When $z$ is not in the divisor of $\Psi(F)$ we
have
\begin{equation}\label{phipsi}
\Phi(z,h;F)=-2\log||\Psi(z,h;F)||^2,
\end{equation}
where $||\cdot||$ is
the Petersson norm, suitably normalized.  Our goal is to evaluate the averages of
$\Phi(F)$ over certain sets of CM points.

\par To define CM points, we consider a rational
splitting
\[
V=V_+\oplus U,
\]
where $V_+$ has signature $(n,0)$ and $U$ has signature $(0,2)$.  This
determines a two-point subset $D_0\subset D$ consisting of $U(\mathbb
R)$ with its two possible orientations.  Let $T=\text{GSpin}(U)$ and $K_T=K\cap
T(\Af)$.  We obtain a zero
cycle
\[
Z(U)_K=T(\qq)\backslash\Big(D_0\times T(\Af)/K_T\Big)\hookrightarrow X_K,
\]
which we regard as a set of CM points inside of $X_K$.  The main
theorem is
\begin{theorem}\label{mainintro}
(i) $\Phi(F)$ is finite at all CM points. \\ \\
(ii) There exist explicit constants $\kappa_\eta(m)$ such that
\begin{equation}\label{introphi}
\sum_{z\in
Z(U)_K}\Phi(z;F)=\frac{4}{\emph{vol}(K_T)}\sum_\eta\sum_{m\geq 0}c_\eta(-m)\kappa_\eta(m),
\end{equation}
where the $c_\eta(-m)$'s are the negative Fourier coefficients of $F$.
\end{theorem}
\noindent
Now using relation (\ref{phipsi}) we obtain
\begin{cor}\label{corPsi}
When $Z(U)_K$ does not meet the divisor of $\Psi(F)$, we have
\begin{equation}\label{intropsi}
\sum_{z\in
Z(U)_K}\log||\Psi(z;F)||^2=\frac{-2}{\emph{vol}(K_T)}\sum_\eta\sum_{m\geq 0}c_\eta(-m)\kappa_\eta(m).
\end{equation}
\end{cor}
\noindent
When $Z(U)_K$ meets the divisor of $\Psi(F)$, it remains to give an
interpretation of Theorem \ref{mainintro} in terms of the function $\Psi(F)$.

\par The constants $\kappa_\eta(m)$ come from Eisenstein series on
$SL_2$.  The quantity in the left hand side of (\ref{introphi}) can be
written as an integral
\[
\int_{\mathbb S(U)}\Phi(z_0,h;F)dh=\int_{\mathbb S(U)}\int_{\Gamma\backslash\h}^{\bullet}((\;F(\tau),\theta(\tau,z_0,h)\;))
v^{-2}dudvdh,
\]
where $\mathbb S(U)=SO(U)(\qq)\backslash SO(U)(\Af)$ and $z_0\in D_0$.  Here we write the theta function as a tensor product
\[
\theta(\tau,z_0,h)=\theta_+(\tau,z_0)\otimes\theta_-(\tau,z_0,h)
\]
of the theta functions for $V_+$ and $U$, respectively.  Then we use the contraction map $\langle\cdot,\theta_+\rangle$ (see section \ref{contract} for details) and write
\[
((\;F(\tau),\theta(\tau,z_0,h)\;))=((\;\langle F,\theta_+\rangle(\tau),\theta_-(\tau,z_0,h)\;)),
\]
where $\langle F,\theta_+\rangle\in S(U(\Af))$.  After some justification, the order of integration (where the inside
integral is regularized) can be switched giving
\begin{equation}\label{introlongint}
\int_{\Gamma\backslash\h}^{\bullet}((\;\langle F,\theta_+\rangle(\tau),\int_{\mathbb S(U)}\theta_-(\tau,z_0,h)dh\;))v^{-2}dudv.
\end{equation}
Then by the Siegel-Weil formula, the integral of $\theta_-(\tau,z_0,h)$
on $\mathbb S(U)$ gives rise
to a \emph{coherent} Eisenstein series, $E(\tau,s;-1)$, of weight $-1$.  For the definition of the term coherent, see \cite{GrossCS}.  With
this Eisenstien series, we can write (\ref{introlongint}) as
\begin{equation}\label{introwithE}
\int_{\Gamma\backslash\h}^{\bullet}((\;\langle F,\theta_+\rangle(\tau),E(\tau,0;-1)\;))v^{-2}dudv.
\end{equation}
Using Maass operators we relate $E(\tau,s;-1)$ to another
Eisenstein series, $E(\tau,s;+1)$, of weight $+1$ via
\[
E(\tau,s;-1)v^{-2}=\frac{-4i}{s}\frac{\partial}{\partial\bar{\tau}}\left
\{E(\tau,s;+1)\right\}.
\] 
One phenomenon that is very specific to the case
of signature $(0,2)$ is that the resulting Eisenstein series
$E(\tau,s;+1)$ is incoherent.  Hence, $E(\tau,s;+1)$ satisfies an odd
functional equation with respect to $s\mapsto -s$, and, therefore,
vanishes at $s=0$.
The integral (\ref{introwithE}) can be evaluated using a Stokes' Theorem argument and some convergence
estimates about the Fourier coefficients of $E(\tau,s;+1)$. This leads
to the constants $\kappa_\eta(m)$ as follows.

\par For $V=V_+\oplus U$ and $L\subset V$, let $L_+=V_+\cap L$ and
$L_-=U\cap L$.  If $\mu\in L_-^\vee/L_-$ and
$\varphi_\mu=\text{char}(\mu+L_-)$ we write
\[
E(\tau,s;\varphi_\mu,+1)=\sum_mA_\mu(s,m,v)\mathbf q^m,
\]
where the Fourier coefficients have Laurent expansions
\[
A_\mu(s,m,v)=b_\mu(m,v)s+O(s^2).
\]
In order to define $\kappa_\eta(m)$, we first define
\[
\kappa_\mu^U(m)=\begin{cases} 
\lim_{v\to\infty}b_\mu(m,v) &\text{if $m>0$,}\\
k_0(0)\varphi_\mu(0) &\text{if $m=0$},\\
0 &\text{if $m<0$},
\end{cases}
\]
where $k_0(0)$ is a constant which depends on the space $U$ (see
Definition \ref{kappadef}).  Let
\[
L^\vee=\bigcup_\eta(\eta+L),\;L=\bigcup_\lambda(\lambda+L_++L_-)
\]
and write $\eta=\eta_++\eta_-, \lambda=\lambda_++\lambda_-$.  Then we
define
\[
\kappa_\eta(m)=\sum_\lambda\sum_{x\in\eta_++\lambda_++L_+}\kappa_{\eta_-+\lambda_-}^U(m-Q(x)).
\]

\par The space $U$ is a rational quadratic space of signature $(0,2)$,
so $U\simeq k$ for an imaginary quadratic field $k$ and the quadratic
form is just a negative multiple of the norm-form.  When $k$ has odd
discriminant and $m\neq 0$, then
\[
\frac{-2}{\text{vol}(K_T)}\kappa_\eta(m)
\]
is the logarithm of an integer.  Thus, if $F$ has
$c_0(0)=0$, so that $\Psi(F)$ is a meromorphic function, then Corollary \ref{corPsi} shows that
\begin{equation}\label{introprod}
\prod_{z\in Z(U)_K}||\Psi(z;F)||^2
\end{equation}
is a rational number.  Moreover, if all of the negative Fourier
coefficients of $F$ are non-negative, then (\ref{introprod}) is
actually an integer.  In the case of signature $(2,2)$, similar
results were obtained in \cite{BrY} for certain rational functions and
CM points on a Hilbert modular surface.  If $c_0(0)\neq 0$, then there is a
transcendental factor
\[
(4\pi d)^{-1}\,e^{2\frac{L'(0,\chi)}{L(0,\chi)}}
\]
appearing in (\ref{introprod}), 
which is related to Shimura's period invariant, \cite{ShimuraAnn},\cite{GrossCS},\cite{Yoshida}, for the CM points in the $0$-cycle $Z(U)_K$.  This factor arises from the trivialization over the CM cycle of the line bundle of which $\Psi(F)$ defines a section.

\par We can say a little bit more about the rational number appearing in
(\ref{introprod}).  The formulas we obtain for $\kappa_\eta(m)$ tell
us the explicit factorization of the rational part of
(\ref{introprod}).  Then, as a consequence of Corollary \ref{corPsi}, we
are able to state a Gross-Zagier type of theorem about which primes
can occur in the factorization.  For $F$ as in (\ref{introFourier}),
define
\[
m_{\text{max}}=\text{max}\{m>0\mid c_\eta(-m)\neq 0\;\text{for some}\;\eta\}.
\] 
\begin{theorem}\label{introsmallprimes}
Let $-d$ be an odd fundamental discriminant and assume $U\simeq
k=\qq(\sqrt{-d})$.  Also assume that $L_-\simeq\a$ for an $\O_k$-ideal $\a$.  Then the only primes which occur in the
factorization of the rational part of
\[
\prod_{z\in Z(U)_K}||\Psi(z;F)||^2
\]
are \\ \\
(i) $q$ such that $q\mid d$, \\ \\
(ii) $p$ inert in $k$ with $p\leq dm_{\emph{max}}$.
\end{theorem}

\par As mentioned in Theorem \ref{mainintro}, one striking phenomenon that occurs in this paper is that the
regularized theta lift $\Phi(F)$ is always finite!  This is interesting since the
Borcherds form $\Psi(F)$ can have zeroes or poles, and
(\ref{phipsi}) only holds when the right hand side is finite.  Considering this, one might
say that the theta lift is \emph{over-regularized}, and it would be
interesting to find the analog of Corollary \ref{corPsi} when $Z(U)_K$
meets the divisor of $\Psi(F)$.

\par There exists lots of recent work on singular moduli, particularly
traces of singular moduli (e.g. \cite{AO}, \cite{BrF} and \cite{Zagier98}).  By considering the case of signature
$(1,2)$, Theorem 1.3 of \cite{GZ} can be recovered
from Theorem \ref{mainintro}.  The appropriate quadratic space is  
\[
V=\{x\in M_2(\zz)\mid\text{tr}(x)=0\}
\]
with $Q(x)=\text{det}(x)$.  For a particular choice of $F$,
\[
\prod_{z\in Z(U)_K}\Psi(z;F)=\prod_{[\tau_1],[\tau_2]}\Big(j(\tau_1)-j(\tau_2)\Big),
\]
where $\tau_1$ and $\tau_2$ are CM points with relatively prime
fundamental discriminants and $[\tau_i]$ denotes an equivalence class
modulo $SL_2(\zz)$.  The right hand side of (\ref{intropsi}) then
gives the same factorization as in \cite{GZ}.  We will discuss this
new proof of Gross-Zagier in a subsequent paper.

\section{Main theorem in the case of signature $(0,2)$}

\subsection{Basic Setup}\label{basicsetup}

We begin by introducing some notation and relevant background
material, and we refer the reader to section 1 of \cite{IoBf} for more details.  Let $V$ be a vector space over $\mathbb Q$ of dimension $n+2$ with
quadratic form $Q$, of signature $(n,2)$, on $V$.  Let $D$ be the space of
oriented negative-definite 2-planes in $V(\mathbb R)$.  For $z\in D$, let
$\text{pr}_z:V(\mathbb R)\rightarrow z$ be the projection map and, for $x\in
V(\mathbb R)$, let $R(x,z)=-(\text{pr}_z(x),\text{pr}_z(x))$.  Then we define
\[
(x,x)_z=(x,x)+2R(x,z),
\]
and our Gaussian for $V$ is the function
\[
\varphi_\infty(x,z)=e^{-\pi(x,x)_z}.
\] 
For $\tau\in\h,\tau=u+iv,$ let 
\[
g_\tau=\left(\begin{matrix}
1&u\\
 &1
\end{matrix}\right)\left(\begin{matrix}
v^{\frac{1}{2}}& \\
 &v^{-\frac{1}{2}}
\end{matrix}\right),
\]
and $g'_\tau=(g_\tau,1)\in Mp_2(\mathbb R)$.
Let $l=\frac{n}{2}-1$, $G=SL_2$ and $\omega$ be the Weil
representation of the metaplectic group $G^\prime_\mathbb A$ on $S(V(\mathbb 
A))$, the
Schwartz space of $V(\mathbb A)$.  If $H=\text{GSpin}(V)$, then for the linear 
action of
$H(\mathbb A_f)$ we write
$\omega(h)\varphi(x)=\varphi(h^{-1}x)$ for $\varphi\in S(V(\mathbb
A_f))$.  If $z\in D$ and $h\in H(\mathbb A_f)$, 
we have the linear functional on $S(V(\mathbb A_f))$ given by
\begin{equation}\label{thetaphi}
\varphi\longmapsto\theta(\tau,z,h;\varphi)=v^{-\frac{l}{2}}\sum_{x\in
V(\mathbb
Q)}\omega(g^\prime_\tau)(\varphi_\infty(\cdot,z)\otimes\omega(h)\varphi)(x).
\end{equation}

\par Let $L\subset V$ be a lattice with dual
\[
L^\vee=\{x\in V\mid (x,L)\subseteq\zz\}
\]
and let $S_L\subset S(V(\mathbb A_f))$
be the space of functions with support
in $\hat{L}^\vee$ and constant on cosets of $\hat{L}$, where $\hat{L}=L\otimes_\zz\hat{\zz}$.  We remark that
$S_L$ is finite dimensional and has a natural basis given by
\[
\{\varphi_\eta=\text{char}(\eta+L)\mid \eta\in L^\vee/L\}.
\]
We also have
\[
S(V(\Af))=\lim_{\substack{\longrightarrow \\ L}}S_L.
\]
\par Let $\Gamma'=Mp_2(\zz)$ be the full inverse image of
$\Gamma=SL_2(\zz)\subset G(\mathbb R)$ in $G'_{\mathbb R}$.  For
$F:\h\to S_L$, the Fourier expansion of $F$ can be written
\begin{equation}\label{FE}
F(\tau)=\sum_\eta
F_\eta(\tau)\varphi_\eta=\sum_\eta\sum_m c_\eta(m)\mathbf q^m\varphi_\eta.
\end{equation}
\begin{df}\label{typeomega}
We say $F:\h\to S_L$ is a 
weakly holomorphic modular form of weight
$1-\frac{n}{2}$ and type $\omega$ for $\Gamma'$ if \\
(i) $F(\gamma'\tau)=(c\tau+d)^{1-\frac{n}{2}}\omega(\gamma')(F(\tau))$
for all $\gamma'\in\Gamma',$ where $\gamma'\mapsto\gamma=\left(\begin{smallmatrix}
a&b \\
c&d
\end{smallmatrix}\right)\in\Gamma$, \\
(ii) $F$ is meromorphic at the cusp, i.e., only a finite number of
the $c_\eta(m)$'s with $m<0$ are non-zero.
\end{df}
\noindent
Note that when $n$ is even, $\omega$ is a representation of
$G_{\mathbb A}$ and we can just work with $\Gamma$.  
The Fourier expansion in (\ref{FE}) is essentially the Fourier expansion given in \cite{Bor98}, where
in that paper he works with group ring elements $\mathbf
e_\eta\in\mathbb C[L^\vee/L]$ instead of the Schwartz functions $\varphi_\eta$.
Since the theta function $\theta(\tau,z,h)$ is a linear functional and
$F(\tau)\in S(V(\Af))$, we can define the $\mathbb C$-bilinear pairing 
\[
((\;F(\tau),\theta(\tau,z,h)\;))=\theta(\tau,z,h;F(\tau)).
\] 
In terms of the Fourier expansion of $F$, this is 
\[
((\;F(\tau),\theta(\tau,z,h)\;))=\sum_\eta F_\eta(\tau)\theta(\tau,z,h;\varphi_\eta).
\] 
Note that as a function of $\tau$, the above pairing is
$\Gamma$-invariant (with a pole at the cusp) since the weights of $\theta$ and $F$ cancel and
their types are dual.
Using this pairing we define
\[
\Phi(z,h;F):=\int_{\Gamma\backslash\h}^{\bullet}((\;F(\tau),\theta(\tau,z,h)\;))
d\mu(\tau),
\]
where $d\mu(\tau)=v^{-2}dudv$ and the integral is regularized as in
\cite{Bor98}.  The regularization is defined
by
\begin{equation}\label{regdef}
\int_{\Gamma\backslash\h}^{\bullet}\phi(\tau)d\mu(\tau)=
\underset{\sigma=0}{\textnormal{CT}}\bigg\{\lim_{t\to\infty}\int_{\mathcal 
F_t}\phi(\tau)v^{-\sigma}d\mu(\tau)\bigg\},
\end{equation}
where we take the constant term in the Laurent expansion at $\sigma=0$
of 
\[
\lim_{t\to\infty}\int_{\mathcal F_t}\phi(\tau)v^{-\sigma}d\mu(\tau),
\]
defined initially for Re$(\sigma)$ sufficiently large.  Here $\mathcal F$ is the
usual fundamental domain for the action of $\Gamma$ on $\h$ and 
\[
\mathcal
F_t=\{\tau\in\mathcal F\mid\text{Im}(\tau)\leq t\}
\]
is the truncated fundamental domain.
  
\subsection{Borcherds Forms}\label{borcherdsforms}

The space $D$ is a bounded symmetric domain.  It can be viewed as an
open subset $\mathcal Q_-$ of a quadric in $\mathbb P(V(\mathbb C))$.
Explicitly,
\[
D\simeq\mathcal Q_-=\{w\in V(\mathbb C)\mid (w,w)=0,
(w,\bar{w})<0\}/\mathbb C^\times,
\]
where the explicit isomorphism is $[z_1,z_2]\mapsto w=z_1+iz_2$ for a
properly oriented basis $[z_1,z_2]$.  Assume $K$ is a compact open
subgroup of $H(\Af)$ such that $H(\mathbb A)=H(\qq)H(\mathbb R)^+K$,
where $H(\mathbb R)^+$ is the identity component of $H(\mathbb R)$.
Define
\[
X_K:=H(\qq)\backslash\big(D\times H(\Af)/K\big).
\]
This is the set of complex points of a quasi-projective variety
rational over $\qq$, and if $\Gamma_K=H(\qq)\cap H(\mathbb R)^+K$,
then $X_K\simeq\Gamma_K\backslash D^+$, where $D^+\subset D$ is the
subset of positively oriented 2-planes.
\par Let $\mathcal L_D$ be the restriction to $D\simeq\mathcal Q_-$ of
the tautological line bundle on $\mathbb P(V(\mathbb C))$.  From this
we get a holomorphic line bundle $\mathcal L$ on $X_K$ equipped with a natural norm, $||\cdot||_{\text{nat}}$, called the Petersson norm.
Assume we have
\[
V(\mathbb R)=V_0+\mathbb Re+\mathbb Rf,
\]
where $e$ and $f$ are such that $(e,f)=1, (e,e)=0=(f,f)$.  Then
sig$(V_0)=(n-1,1)$ and for the negative cone
\[
\mathcal C=\{y\in V_0\mid (y,y)<0\},
\]
we have
\[
D\simeq\mathbb D:=\{z\in V_0(\mathbb C)\mid y=\text{Im}(z)\in\mathcal C\}.
\]
The explicit isomorphism is 
\[
\mathbb D\to V(\mathbb C),\;z\mapsto w(z):=z+e-Q(z)f
\]
composed with projection to $\mathcal Q_-$.  The map $z\mapsto w(z)$
can be viewed as a holomorphic section of $\mathcal L_D$.
\par We now define the notion of a modular form on $D\times H(\Af)$.
\begin{df}
A modular form on $D\times H(\Af)$ of weight $m\in\frac{1}{2}\zz$ is a function
$\Psi:D\times H(\Af)\to \mathbb C$ such that
\begin{enumerate}
\item $\Psi(z,hk)=\Psi(z,h)$ for all $k\in K$,
\item $\Psi(\gamma z,\gamma h)=j(\gamma,z)^m\Psi(z,h)$ for all
$\gamma\in H(\qq)$, where $j(\gamma,z)$ is an automorphy factor.
\end{enumerate}
\end{df}
\noindent
Meromorphic modular forms on $D\times H(\Af)$ of weight $m\in\zz$ can be identified with
meromorphic sections of $\mathcal L^{\otimes m}$.  If $\Psi$ is such a
meromorphic modular form, then the Petersson norm of the section
$(z,h)\mapsto \Psi(z,h)w(z)^{\otimes m}$ associated to $\Psi$ is
\[
||\Psi(z,h)||_{\text{nat}}^2=|\Psi(z,h)|^2|y|^{2m}.
\]
For reasons we will see below, we renormalize $||\cdot||_{\text{nat}}$
and instead work with the following norm
\[
||\Psi(z,h)||^2:=||\Psi(z,h)||^2_{\text{nat}}\left(2\pi e^{\Gamma'(1)}\right)^{\!m}.
\]
The ``extra" constant in the metric here is related to that occuring in \cite{Faltht}.  Borcherds proved that the regularized integral $\Phi(z,h;F)$ satisfies
the equation
\begin{align}\label{phiandpsi}
\Phi(z,h;F)&=-2\log||\Psi(z,h;F)||_{\text{nat}}^2-c_0(0)(\log(2\pi)+\Gamma'(1))\\
&=-2\log||\Psi(z,h;F)||^2\notag
\end{align}
for a meromorphic modular form $\Psi(F)$ on $D\times H(\Af)$ of weight
$m=\frac{1}{2}c_0(0)$ when $z$ does not lie in the divisor of
$\Psi(F)$.
\begin{remark}
In fact, $\Phi(F)$ may still be finite for $z\in D$ even if $z$
lies in the divisor of $\Psi$.  This value of $\Phi(F)$ must have another
meaning there.
\end{remark}
\begin{df}
A Borcherds form $\Psi(F)$ is a meromorphic modular form
on $D\times H(\Af)$ which arises (via (\ref{phiandpsi})) from the
regularized theta lift of a modular form $F$.
\end{df}

\subsection{CM Points}\label{CMsection}

Assume that
we have a rational splitting
\[
V=V_+\oplus U,
\]
where $V_+$ has signature $(n,0)$ and $U$ has signature $(0,2)$.  This
determines a two-point subset $\{z_0^\pm\}=D_0\subset D$ given by $U(\mathbb
R)$ with its two orientations.  For $z_0\in D_0$, we are interested in 
computing the
integral
\begin{equation}\label{wantcompute}
\int_{SO(U)(\qq)\backslash SO(U)(\Af)}\Phi(z_0,h;F)dh.
\end{equation} 
Let $T=\text{GSpin}(U)$ and note there is a natural homomorphism $T\to H$.
Let $K$ be as in section
\ref{borcherdsforms} and define
$K_T=K\cap T(\Af)$.  Consider the set of CM
points
\[
Z(U)_K:=T(\qq)\backslash\Big(D_0\times T(\Af)/K_T\Big)\hookrightarrow X_K.
\]
We want to compute
\[
\text{vol}(K_T)\sum_{z\in Z(U)_K}\Phi(z;F)=-2\int_{SO(U)(\qq)\backslash SO(U)(\Af)}\Phi(z_0,h;F)dh.
\]
Note that after normalizing by the volume of $K_T$, this expression is
independent of the choice of $K$.

\subsection{Convergence Questions and Regularization}

First we consider the case when $n=0$ and our space $V=U$ is
negative-definite.  In this case, $D=D_0$, the Gaussian is
$\varphi_\infty(x)=e^{\pi(x,x)}$ and the theta function is
\begin{equation}\label{zerotheta}
\theta(\tau,z_0,h;\varphi)=v^{\frac{1}{2}}\sum_{x\in U(\mathbb
Q)}\omega(g_\tau^\prime)e^{\pi(x,x)}\varphi(h^{-1}x),
\end{equation}
for any $\varphi\in S(U(\Af))$.  When $n=0$ and we have a lattice $L\subset
U$ we write $\mu\in L^\vee/L$ and $\varphi_\mu=\text{char}(\mu+L)$.
Let $F(\tau)$ be a weakly holomorphic modular form of
weight 1 valued in $S_L$, and let
\begin{equation}\label{FexpansionL}
F(\tau)=\sum_\mu
F_\mu(\tau)\varphi_\mu=\sum_\mu\sum_{m\in\mathbb
Q}c_\mu(m)\mathbf q^m\varphi_\mu,
\end{equation}
where $\mu$ runs over $L^\vee/L$.  We assume $c_\mu(m)\in\zz$ for $m\leq 0$.  The
functions $F_\mu$ are meromorphic modular forms with some real multiplier for
a congruence subgroup of $SL_2(\zz)$, and it
will be very useful to know how large their Fourier coefficients can
be.  
\begin{lem}\label{wt1bound}
Assume $m_\mu\in\zz$ is such that $c_\mu(m_\mu)\neq 0$ and
$c_\mu(m)=0$ for all $m<
m_\mu$.  Then there are constants $C$ and $C'$ such that, for $m>0$,
\[
\left|c_\mu(m)\right|\leq C'\left((-m_\mu+2)(m-m_\mu)^6+m^6e^{C\sqrt{m}}\right),
\]
where $C$ depends on $m_\mu$ and on the
multiplier and $C'$ depends on the polar part of $F_\mu$.
\end{lem}
\begin{proof}
The cusp form of weight $12$, $(2\pi)^{-12}\Delta(\tau)=\mathbf
q\prod_{n=1}^\infty(1-\mathbf q^n)^{24}$, has Fourier expansion
\[
(2\pi)^{-12}\Delta(\tau)=\sum_{N=1}^\infty\tau(N)\mathbf q^N,
\]
where $|\tau(N)|\leq C_1N^6$ for some constant $C_1$.  Let
$\tilde{\Delta}(\tau)=(2\pi)^{-12}\Delta(\tau)$.  We can look at
$F_\mu/\tilde{\Delta}$, which has weight $-11=1-\frac{24}{2}$.  If
\[
F_\mu/\tilde{\Delta}=\sum_{m=m_\mu-1}^\infty a_\mu(m)\mathbf q^m,
\]
then for $m>0$, (3.38) of \cite{IoBf} tells us there are constants $C_2$
and $C$
such that
\[
|a_\mu(m)|\leq C_2m^{-\frac{25}{4}}e^{C\sqrt{m}},
\]
where $C$ depends on $m_\mu$ and on the multiplier.
We have
\begin{align}\notag
F_\mu(\tau)&=\Big(\sum_{N=1}^\infty\tau(N)\mathbf
q^N\Big)\Big(\sum_{m=m_\mu-1}^\infty a_\mu(m)\mathbf
q^m\Big)\\
&=\sum_{N=1}^\infty\;\sum_{m=m_\mu-1}^\infty
\tau(N)a_\mu(m)\mathbf q^{N+m}\notag\\
&=\sum_{m=m_\mu}^\infty\Bigg[\sum_{N=1}^{m-m_\mu+1}\tau(N)a_\mu(m-N)\Bigg]\mathbf q^m.\notag
\end{align}
Then
\begin{align}\notag
|c_\mu(m)|&=\left|\sum_{N=1}^{m-m_\mu+1}\tau(N)a_\mu(m-N)\right|\\
&=\left|\sum_{N\geq
m}\tau(N)a_\mu(m-N)+\sum_{0<N<m}\tau(N)a_\mu(m-N)\right|\notag\\
&\leq C_1\sum_{N=m}^{m-m_\mu+1}N^6|a_\mu(m-N)|+C_1C_2\sum_{0<N<m}N^6(m-N)^{-\frac{25}{4}}e^{C\sqrt{m-N}}.\notag
\end{align}
We know there is a constant $C_3$
such that $|a_\mu(m)|\leq C_3$ for
$m\in\{m_\mu,\ldots,0\}$, and thus
\begin{align}\notag
|c_\mu(m)|&\leq
 C_1C_3(-m_\mu+2)(m-m_\mu)^6+C_1C_2m^6e^{C\sqrt{m}}\\
&\leq C'\left((-m_\mu+2)(m-m_\mu)^6+m^6e^{C\sqrt{m}}\right),\notag
\end{align}
for some constant $C'$. 
\end{proof}
In the $n=0$ case, the following over-regularization phenomenon occurs:
\begin{prop}\label{alwaysfinite}For $h\in H(\mathbb A_f)$,
\[
\Phi(z_0,h;F)=\int_{\Gamma\backslash\h}^{\bullet}((\;F(\tau),\theta(\tau,z_0,h)\;))d\mu(\tau)
\]
is always finite.
\end{prop}
\begin{proof}
This case corresponds to signature $(2,0)$ in \cite{Bor98}.  In Theorem 6.2 of
\cite{Bor98}, Borcherds points out that $\Phi$ is
nonsingular except along a locally finite set of codimension 2
sub-Grassmannians $\lambda^\perp$, for some negative norm vectors
$\lambda\in L$.  No such vectors exist in signature $(2,0)$.  For
ease of the reader, we
give the proof in our notation.  We have
\begin{equation}\label{defwithF}
\int_{\Gamma\backslash\h}^{\bullet}((\;F(\tau),\theta(\tau,z_0,h)\;))d\mu(\tau)=
\underset{\sigma=0}{\textnormal{CT}}\bigg\{\lim_{t\to\infty}\int_{\mathcal 
F_t}\theta(\tau,z_0,h;F)v^{-\sigma}d\mu(\tau)\bigg\},
\end{equation}
and we can write the
integral on the right hand side of (\ref{defwithF}) as
\[
\int\limits_{1}^{t}\int\limits_{-\frac{1}{2}}^{\frac{1}{2}}\theta(\tau,z_0,h;F)v
^{-\sigma}d\mu(\tau)
+ \int_{\mathcal F_1}\theta(\tau,z_0,h;F)v^{-\sigma}d\mu(\tau).
\]
The integral over the compact set $\mathcal F_1$ is finite and
independent of $t$, so we
just look at the first part.   By \cite{Shintani}, we have
\[
\omega(g'_\tau)e^{\pi(x,x)}=v^{\frac{1}{2}}e(uQ(x))e^{2\pi vQ(x)}, 
\]
where $e(y)=e^{2\pi iy}$.
Then (\ref{zerotheta}) is
\[
\theta(\tau,z_0,h;\varphi)=v\sum_{x\in
U(\mathbb Q)}e(uQ(x))e^{2\pi vQ(x)}\varphi(h^{-1}x),
\]
and so the integral over $\mathcal F_t-\mathcal F_1$ is
\begin{equation}\label{ftminusf1}
\sum_{\mu}\sum_{m\in\mathbb Q}\sum_{x\in U(\mathbb
Q)}c_\mu(m)\varphi_\mu(h^{-1}x)\int\limits_{1}^{t}\int\limits_{-\frac{1}{2}}^{\frac{1}{2}} e(um)e(uQ(x))e^{-2\pi
vm}e^{2\pi vQ(x)}v^{-\sigma-1}dudv.
\end{equation}
\begin{lem}
If $m+Q(x)\notin\zz$, then $c_\mu(m)=0$.
\end{lem}
\begin{proof}
When we consider the
transformation law for $F$, we have
$F(\tau+1)=\omega(T)(F(\tau))$.  That is, for any $x\in
U(\mathbb A_f)$,
\begin{align}\notag
\sum_\mu\sum_mc_\mu(m)\mathbf 
q^me(m)\varphi_\mu(x)&=\omega(T)\left(\sum_\mu\sum_mc_\mu(m)
\mathbf q^m\varphi_\mu\right)(x)\\
&=\sum_\mu\sum_mc_\mu(m)\mathbf 
q^m\omega(T)(\varphi_\mu)(x)\notag\\
&=\sum_\mu\sum_mc_\mu(m)\mathbf q^me(-Q(x))\varphi_\mu(x).\notag
\end{align}
We see
$m+Q(x)\notin\mathbb Z$ implies $c_\mu(m)=0$.  
\end{proof}
\noindent For
$m+Q(x)\in\mathbb Z$,
\[
\int\limits_{-\frac{1}{2}}^{\frac{1}{2}} e(um)e(uQ(x))du=\begin{cases} 
1 &\text{if $m+Q(x)=0,$}\\
0 &\text{otherwise}.
\end{cases}.
\]
Integrating with respect to $u$ in (\ref{ftminusf1}) and 
letting $t\to\infty$ gives
\begin{equation}\label{ugone}
\sum_{\mu}\sum_{\substack{m\in\mathbb Q\\ m\geq 0}}\sum_{\substack{x\in 
U(\mathbb
Q) \\ Q(x)+m=0}}c_\mu(m)\varphi_\mu(h^{-1}x)\int\limits_{1}^{\infty} e^{-4\pi 
mv}v^{-\sigma-1}dv.
\end{equation}
We have $m\geq 0$ since $Q(x)\leq 0$.  When $m=0$, we get
\[
\sum_\mu c_\mu(0)\varphi_\mu(0)\int\limits_{1}^{t} 
v^{-\sigma-1}dv=c_0(0)\frac{1}{\sigma}(1-t^{-\sigma}), 
\]
which equals zero when we take the limit as $t\to\infty$
followed by the constant term at $\sigma=0$.  For $m>0$, (3.35) of
\cite{IoBf} says
\[
\int\limits_{0}^{\infty}e^{-4\pi mv}v^{-\sigma-1}dv\leq
C(\epsilon,\sigma)e^{-4\pi m}
\]
for any $\epsilon$ with $0<\epsilon<4\pi m,$ where the constant
$C(\epsilon,\sigma)$ is uniform in any $\sigma$-halfplane and
independent of $m$.  Using this in (\ref{ugone}), we have
\[
C(\epsilon,\sigma)\sum_\mu\sum_{m>0}c_\mu(m)e^{-4\pi 
m}\sum_{\substack{x\in U(\mathbb
Q) \\ Q(x)+m=0}}\varphi_\mu(h^{-1}x),
\]
which is finite by Lemma \ref{wt1bound}.
\end{proof}

\subsection{Eisenstein Series}\label{subsectioneseries}

Here we give the basic definition of an Eisenstein series and some
related theory when $V$ has signature $(n,2)$ for $n$ even.  What
follows is a summary of the explanations given in \cite{IoBf} for $n$ even,
and we refer the reader to that paper for the more general theory.  Inside of $G_{\mathbb A}$, we
have the subgroups
\[
N_{\mathbb A}=\{n(b)\mid b\in\mathbb A\},\;n(b)=\left(\begin{matrix}
1&b\\
 &1
\end{matrix}\right),
\]
and
\[
M_{\mathbb A}=\{m(a)\mid a\in\mathbb A^\times\},\;m(a)=\left(\begin{matrix}
a& \\
 &a^{-1}
\end{matrix}\right).
\]
Define the quadratic character $\chi=\chi_V$ of $\mathbb
A^\times/\qq^\times$ by
\[
\chi(x)=(x,-\text{det}(V)),
\]
where $\text{det}(V)\in\qq^\times/(\qq^\times)^2$ is the determinant
of the matrix for the quadratic form $Q$ on $V$.  For $s\in\mathbb C$,
let $I(s,\chi)$ be the principal series representation of $G_{\mathbb
A}$.  This space consists of smooth functions $\Phi(s)$ on $G_{\mathbb
A}$ such that
\[
\Phi(n(b)m(a)g,s)=\chi(a)|a|^{s+1}\Phi(g,s).
\]
We have a $G_{\mathbb A}$-intertwining map
\begin{equation}\label{lamsubv}
\lambda=\lambda_V:S(V(\mathbb A))\to I\left(\frac{n}{2},\chi\right),
\end{equation}
where $\lambda(\varphi)(g)=(\omega(g)\varphi)(0)$.  If
$K_\infty=SO(2)$ and $K_f=SL_2(\hat{\zz})$, then a section
$\Phi(s)\in I(s,\chi)$ is called standard if its restriction to
$K_\infty K_f$ is independent of $s$.  The function $\lambda(\varphi)$
has a unique extension to a standard section $\Phi(s)\in I(s,\chi)$
such that $\Phi\left(\frac{n}{2}\right)=\lambda(\varphi)$.  We let $P=MN$ and define the
Eisenstein series associated to $\Phi(s)$ by
\[
E(g,s;\Phi)=\sum_{\gamma\in P_\qq\backslash G_\qq}\Phi(\gamma g,s),
\]
where $G_\qq$ is identified with its image in $G_{\mathbb A}$.  This
series converges for $\text{Re}(s)>1$ and has a meromorphic analytic
continuation to the whole $s$-plane.
\par One step in proving the $(0,2)$-Theorem is to apply Maass operators
to obtain a relation between two Eisenstein series.  Let
\[
X_{\pm}=\frac{1}{2}\left(\begin{matrix}
1&\pm i\\
\pm i&-1
\end{matrix}\right)\in\mathfrak s\mathfrak l_2(\mathbb C).
\]
For $r\in\zz$, let $\chi_r$ be the character of $K_\infty$ defined by
\[
\chi_r(k_\theta)=e^{ir\theta},\;k_\theta=\left(\begin{matrix}
\text{cos}\;\theta&\text{sin}\;\theta\\
-\text{sin}\;\theta&\text{cos}\;\theta
\end{matrix}\right)\in K_\infty.
\]
Let $\phi:G_{\mathbb R}\to\mathbb C$ be a smooth function of weight
$l$, meaning $\phi(gk_\theta)=\chi_l(k_\theta)\phi(g)$, and let
$\xi(\tau)=v^{-\frac{l}{2}}\phi(g_\tau)$ be the corresponding
function on $\h$.  Then $X_\pm\phi$ has weight $l\pm 2$, and the
corresponding function on $\h$ is
\[
v^{-\frac{l\pm 2}{2}}X_\pm\phi(g_\tau)=\begin{cases} 
\left(2i\frac{\partial\xi}{\partial\tau}+\frac{l}{v}\xi\right)(\tau)
&\text{for $+$,}\\
-2iv^2\frac{\partial\xi}{\partial\bar{\tau}}(\tau) &\text{for $-$.}
\end{cases}
\]
\begin{lem}[Lemma 2.7 of \cite{IoBf}]\label{ssk2.7}
Let $\Phi^r_\infty(s)\in I_\infty(s,\chi)$ be the normalized
eigenvector of weight $r$ for the action of $K_\infty$.  Then
\[
X_\pm\Phi^r_\infty(s)=\frac{1}{2}(s+1\pm r)\Phi^{r\pm 2}_\infty(s).
\]
\end{lem}
\noindent
For $\varphi\in S(V(\Af))$, let $E(g,s;\Phi^r_\infty\otimes
\lambda(\varphi))$ be the Eisenstein series of weight $r$ on
$G_{\mathbb A}$ associated to $\varphi$.  For the Gaussian,
$\varphi_\infty(x,z)$, we have
$\lambda(\varphi_\infty)=\Phi^l_\infty\left(\frac{n}{2}\right)$, where $l=\frac{n}{2}-1$.
This means that
\[
X_-E(g,s;\Phi^{l+2}_\infty\otimes\lambda(\varphi))=\frac{1}{2}(s-l-1)E(g,s;\Phi^l_\infty\otimes\lambda(\varphi)).
\]
On $\h$, this translates to
\begin{equation}\label{maassonh}
-2iv^2\frac{\partial}{\partial\bar{\tau}}\Big\{E(\tau,s;\varphi,l+2)\Big\}=\frac{1}{2}\left(s-\frac{n}{2}\right)E(\tau,s;\varphi,l),
\end{equation}
where we write
$E(\tau,s;\varphi,l)=v^{-\frac{l}{2}}E(g_\tau,s;\Phi^{l}_\infty\otimes\lambda(\varphi))$.
One main result we need is the Siegel-Weil formula.
\begin{theorem}[Siegel-Weil formula]\label{thesw}
Let $V$ be a vector space of signature $(n,2)$.  Assume $V$ is
anisotropic or that $\emph{dim}(V)-r_0>2$, where $r_0$ is the Witt
index of $V$.  Then $E(g,s;\varphi)$ is holomorphic at $s=\frac{n}{2}$
and
\[
E\left(g,\frac{n}{2};\varphi\right)=\frac{\alpha}{2}\int_{SO(V)(\qq)\backslash SO(V)(\mathbb
A)}\theta(g,h;\varphi)dh,
\]
where $dh$ is Tamagawa measure on $SO(V(\mathbb A))$, and $\alpha$ is
$2$ if $n=0$ and is $1$ otherwise.
\end{theorem}
\noindent
Here $\theta(g,h;\varphi)$ is defined as in (\ref{thetaphi}) without
$v^{-\frac{l}{2}}$ and with $g$ replacing $g'_\tau$.  The integration
for $SO(U)(\mathbb R)$ is with
respect to the action $h_\infty^{-1}x$ in the argument of
$\varphi_\infty.$  The cases which are omitted in the Siegel-Weil
formula are when $n=1=r_0$ ($V$ is isotropic) and $n=2=r_0$ ($V$ is split).
\par Let us now consider the situation $V=U$, $\text{sig}(U)=(0,2)$.  The
representation we are interested in is $I(0,\chi)$.  This global
principal series is a restricted tensor product of local ones,
\[
I(0,\chi)=\otimes'_v I_v(0,\chi_v).
\]
For the local space $U_v=U(\qq_v)$, define the quadratic character
$\chi_v$ of $\qq_v^\times$ by
\[
\chi_v(x)=(x,-\text{det}(U_v))_v.
\]
Let $R_v(U)$ be the maximal quotient of $S(U_v)$ on which $O(U_v)$
acts trivially.  The following proposition is a special case of
Proposition 1.1 of \cite{SSKAnnals}.
\begin{prop}
(i) If $v\neq\infty$, then 
\[
I_v(0,\chi_v)=R_v(U^+)\oplus R_v(U^-),
\]
where $U^\pm$ has Hasse invariant $\epsilon_v(U^\pm)=\pm 1$. \\
(ii) If $v=\infty$, then
\[
I_\infty(0,\chi_\infty)=R_\infty(U(0,2))\oplus R_\infty(U(2,0)),
\] 
and the spaces $U(0,2)$ and $U(2,0)$ have opposite Hasse invariants.
\end{prop}
\noindent
Recall the notion of an incoherent collection.
\begin{df}
An incoherent collection $\mathcal C=\{\mathcal C_v\}$ of
quadratic spaces is a set of quadratic spaces $\mathcal C_v$ such that
\begin{enumerate}
\item For all $v$, $\emph{dim}_{\qq_v}(\mathcal C_v)=2$, and
$\chi_{\mathcal C_v}=\chi$.
\item For almost all $v$, $\mathcal C_v\simeq U_v$.
\item (Incoherence condition) The product formula fails for the Hasse
invariants:
\[
\prod_v\epsilon_v(\mathcal C_v)=-1.
\]
\end{enumerate}
\end{df}
\noindent
Then we have, cf. (2.10) in \cite{SSKAnnals},
\[
I(0,\chi)\simeq\left(\bigoplus_{U'}\Pi(U')\right)\oplus\left(\bigoplus_{\mathcal
C}\Pi(\mathcal C)\right)
\]
as a sum of two irreducible pieces defined as follows.  $U'$ runs over
all global quadratic spaces of dimension $2$ with $\chi_{U'}=\chi$,
while $\mathcal C$ runs over all
incoherent collections of dimension $2$ and character $\chi$, and
\[
\Pi(U')=\otimes'_v R_v(U'),\;\Pi(\mathcal C)=\otimes'_v R_v(\mathcal C).
\]
For $\lambda=\lambda_U$ as in
(\ref{lamsubv}), we have
$\lambda(\varphi_\infty)=\Phi^{-1}_\infty(0)$, where
$\Phi^{-1}_\infty$ is the normalized eigenvector of weight $-1$ for
the action of $K_\infty$.  From the theory of principal series
representations, we have
$\Phi^{-1}_\infty(0)\in R_\infty(U(0,2))$ and
$\Phi^{1}_\infty(0)\in R_\infty(U(2,0))$.  Then Lemma \ref{ssk2.7}
implies
\begin{equation}\label{shift}
X_+\Phi^{-1}_\infty(s)=\frac{1}{2}s\Phi^{1}_\infty(s),
\end{equation}
so we see that the Maass operator $X_+$ shifts the coherent
Eisenstein series\\ $E(g,s;\Phi^{-1}_\infty\otimes\lambda(\varphi))$ to
the \emph{incoherent} Eisenstein series
$E(g,s;\Phi^{1}_\infty\otimes\lambda(\varphi))$.  Theorem 2.2 of
\cite{SSKAnnals} then tells us that
\[
E(g,0;\Phi^{1}_\infty\otimes\lambda(\varphi))=0.
\]

\subsection{The $(0,2)$-Theorem}\label{subsectionn=0}

The integral we want to
compute is
\begin{equation}\label{mainintegral}
\int_{SO(U)(\qq)\backslash SO(U)(\Af)}\Phi(z_0,h;F)dh,
\end{equation}
which is equal to
\begin{equation}\label{insidedot}
\int_{SO(U)(\qq)\backslash 
SO(U)(\Af)}\int_{\Gamma\backslash\h}^{\bullet}((\;F(\tau),\theta(\tau,z_0,h)\;))
d\mu(\tau)dh.
\end{equation}
As in \cite{IoBf}, we would like to be able to switch the order of
integration, where the inside integral is regularized.  That is, we
want (\ref{insidedot}) to equal
\[
\int_{\Gamma\backslash\h}^{\bullet}((\;F(\tau),\int_{SO(U)(\qq)\backslash 
SO(U)(\Af)}\theta(\tau,z_0,h)dh\;))d\mu(\tau).
\]
Note that $F:\h\to S_L$ implies $F(\tau)\in S(U(\mathbb A_f))^{K}$,
where 
\[
K=\{h\in H(\Af)\mid
h(\lambda+L)=\lambda+L,\forall\lambda\in L^{\vee}/L\}
\]
is an open subset of $H(\Af)$.  

\par Before we justify the
interchange of integrals, we need to make some remarks about our specific
case.  For a reference on Clifford algebras, see \cite{Chevalley} or \cite{Hermann}.  The 
Clifford algebra
$C(U)$ can be written as $C(U)=C^0(U)\oplus C^1(U)$, where
$C^0(U)$ and $C^1(U)$ are the even and odd parts, respectively.
$C^0(U)^\times$ acts on $C^1(U)$ by conjugation.  Assume
$U$ has basis $\{u,v\}$ with $Q(u)=a, Q(v)=b$ and $(u,v)=0$. Then $C(U)$ is
spanned by $\{1,u,v,uv\}$ with $C^0(U)=\text{span}\{1,uv\}$ and
$C^1(U)=\text{span}\{u,v\}$.  By definition,
\[
H=\{g\in C^0(U)^\times\mid gUg^{-1}=U\}.
\]
Since $C^1(U)=U$, $H=C^0(U)^\times$.  In $C^0(U)$ we have
$(uv)^2=-ab$, so if $k=\mathbb Q\left(\sqrt{-ab}\right)$, then $H\simeq 
k^\times$.
This means $SO(U)\simeq k^1$ and $k^\times\to k^1$ is the map
\[
x\mapsto\frac{x}{x^\sigma}
\]
by Hilbert's Theorem 90.  We have the
exact sequence
\[
1\to Z\to H\to SO(U)\to 1,
\]
where $H(\mathbb A_f)\simeq k^\times_{\mathbb A_f}, H(\mathbb
Q)\simeq k^\times, Z(\mathbb A_f)\simeq\mathbb Q^\times_{\mathbb
A_f}$ and $Z(\mathbb Q)\simeq\mathbb Q^\times.$  
\begin{lem}\label{Uisk}
For any negative-definite space $U$ with quadratic form $Q$ of
signature $(0,2)$, we realize $U\simeq k$ for an imaginary quadratic
field $k$ and $Q$ is given by a negative multiple of the norm-form.
\end{lem}
\noindent
If $B(h)$ is a function on
$H(\Af)$ which only depends on the image of $h$ in $SO(U)(\Af)$, then
we can view $B$ as a function on $SO(U)(\Af)$ as well.    
\begin{lem}\label{integraltosum}
Let $B(h)$ be a function on
$H(\Af)$ depending only on the image of $h$ in $SO(U)(\Af)$.
Assume $B$ is invariant under $K$ and $H(\qq)$.  Then 
\[
\int_{SO(U)(\qq)\backslash SO(U)(\mathbb A_f)}B(h)dh=\emph{vol}(K)\sum_{h\in 
H(\qq)\backslash H(\Af)/K}B(h),
\]
and the sum is finite. 
\end{lem}
\begin{proof}
We have the exact sequence
\[
1\to k^0_{\mathbb A}\to k^\times_{\mathbb A}\to\mathbb R^\times_+\to 1,
\]
where the map to $\mathbb R^\times_+$ is the absolute value map.  By
the product formula, $k^\times\subset k^0_{\mathbb A}$ and we know
$k^\times\backslash k^0_{\mathbb A}$ is compact.
\begin{lem}\label{kcrossa}
$k^\times_{\mathbb A}=\qq^\times_{\mathbb A}k^0_{\mathbb A}$.
\end{lem}
\begin{proof}
$\qq^\times_{\mathbb A}$ injects into $k^\times_{\mathbb A}$ and also
maps onto $\mathbb R^\times_+$.  So if $(a)\in k^\times_{\mathbb A}$
then $\exists (b)\in\qq^\times_{\mathbb A}$ with $|(b)|=|(a)|$.  Then
$(b)\in\qq^\times_{\mathbb A}\subset k^\times_{\mathbb A}$ implies
$k^0_{\mathbb A}(b)=k^0_{\mathbb A}(a)$, so $(a)\in\qq^\times_{\mathbb 
A}k^0_{\mathbb A}$.
\end{proof}
\noindent
Lemma \ref{kcrossa} implies
\[
k^\times\backslash k^0_{\mathbb A}\twoheadrightarrow k^\times\mathbb
Q^\times_{\mathbb A}\backslash k^\times_{\mathbb A},
\]
and so $k^\times\mathbb
Q^\times_{\mathbb A}\backslash k^\times_{\mathbb A}$ is also compact.
The set we integrate over is
\[
SO(U)(\qq)\backslash SO(U)(\mathbb A_f)=H(\qq)\backslash
H(\Af)/Z(\Af)\simeq k^\times\qq^\times_{\Af}\backslash k^\times_{\Af}.
\]
This is compact since $k^\times\mathbb
Q^\times_{\mathbb A}\backslash k^\times_{\mathbb A}$ maps onto it.
Then $K$ is open and $K\supset Z(\Af)$ so $H(\qq)\backslash
H(\Af)/K$ is finite.  The volume term appears since $B$ is $K$-invariant.
\end{proof}
\begin{prop}\label{switchorder}
\begin{align}\notag
&\int_{SO(U)(\qq)\backslash SO(U)(\mathbb
A_f)}\int_{\Gamma\backslash\h}^{\bullet}((\;F(\tau),\theta(\tau,z_0,h)\;))d\mu(\tau)dh\\
&=\int_{\Gamma\backslash\h}^{\bullet}((\;F(\tau),\int_{SO(U)(\mathbb
Q)\backslash SO(U)(\Af)}\theta(\tau,z_0,h)dh\;))d\mu(\tau).\notag
\end{align}
\end{prop}
\begin{proof}
The main point is that since $F(\tau)\in S(U(\mathbb A_f))^{K}$, we know
\[
\int_{\Gamma\backslash\h}^{\bullet}((\;F(\tau),\theta(\tau,z_0,h)\;))d\mu(\tau)
\]
is $K$-invariant.  So if we let 
\[
B(h)=\int_{\Gamma\backslash\h}^{\bullet}((\;F(\tau),\theta(\tau,z_0,h)\;))d\mu(\tau),
\]
then Lemma \ref{integraltosum} says
\begin{align}\notag
\int_{SO(U)(\qq)\backslash SO(U)(\mathbb
A_f)}B(h)dh&=\text{vol}(K)\sum_{h\in H(\qq)\backslash
H(\Af)/K}B(h)\\
&=\int_{\Gamma\backslash\h}^{\bullet}\text{vol}(K)\sum_{h\in H(\qq)\backslash
H(\Af)/K}\theta(\tau,z_0,h;F(\tau))d\mu(\tau),\label{suminside}
\end{align}
since the sum is finite.  Now apply Lemma \ref{integraltosum} again to 
$\theta(\tau,z_0,h;F(\tau))$
and (\ref{suminside}) is
\[
=\int_{\Gamma\backslash\h}^{\bullet}((\;F(\tau),\int_{SO(U)(\mathbb
Q)\backslash SO(U)(\Af)}\theta(\tau,z_0,h)dh\;))d\mu(\tau).
\]
\end{proof}
\noindent 
The quadratic space $U$ is anisotropic, so we can apply Theorem \ref{thesw}.  This tells us
that for any $\varphi\in S(U(\mathbb A))$, 
\begin{equation}\label{intequalsE}
\int_{SO(U)(\qq)\backslash SO(U)(\mathbb 
A)}\theta(\tau,z_0,h;\varphi)dh=v^{\frac{1}{2}}E(g_\tau,0;\varphi,-1),
\end{equation}
where $E(g_\tau,s;\varphi,-1)$ is a coherent Eisenstein series of
weight $-1$.  Since $\theta(\tau,z_0,h)$ is $SO(U)(\mathbb
R)$-invariant, it suffices to integrate over $SO(U)(\qq)\backslash
SO(U)(\Af)$.  We choose a factorization for the measure
$dh=dh_\infty\times dh_f$ such that $\text{vol}(SO(U)(\mathbb R))=1$. 
\begin{lem}\label{EwithAf}
\begin{align}\notag
&(i)\quad \int_{SO(U)(\qq)\backslash SO(U)(\Af)}\theta(\tau,z_0,
h_f;\varphi)dh_f=v^{\frac{1}{2}}E(g_\tau,0;\varphi,-1).\\
&(ii)\quad \emph{vol}(K)^{-1}=\frac{1}{2}(\#(H(\qq)\backslash 
H(\Af)/K)).\notag
\end{align}
\end{lem}
\noindent
We let 
\[
E(\tau,s;-1):=v^{\frac{1}{2}}E(g_\tau,s;-1).
\]
Then for (\ref{mainintegral}) we have
\begin{equation}\label{mainintEis}
\int_{SO(U)(\qq)\backslash SO(U)(\mathbb
A_f)}\Phi(z_0,h;F)dh=\int_{\Gamma\backslash\h}^{\bullet}((\;
F(\tau),E(\tau,0;-1)\;))d\mu(\tau).
\end{equation}
For $F$ as in (\ref{FexpansionL}), the right hand side of
(\ref{mainintEis}) is
\begin{equation}\label{Eonphi}
\int_{\Gamma\backslash\h}^{\bullet}((\;
F(\tau),E(\tau,0;-1)\;))d\mu(\tau)=\int_{\Gamma\backslash\h}^{\bullet}\sum_\mu 
F_\mu(\tau)E(\tau,0;\varphi_\mu,-1)v^{-2}dudv.
\end{equation}
Let
\[
I(s,t):=\int_{\mathcal F_t}\sum_\mu 
F_\mu(\tau)E(\tau,s;\varphi_\mu,-1)v^{-2}dudv.
\]
In order to state the main theorem of this chapter, we view
$U\simeq k=\qq(\sqrt{-d})$, where $-d$ is the discriminant of $k$, and
let $\chi_d$ be the character of $\qq^\times_{\mathbb A}$ defined by
$\chi_d(x)=(x,-d)_{\mathbb A}$.  We define the normalized $L$-series
\[
\Lambda(s,\chi_d)=\pi^{-\frac{s+1}{2}}\Gamma\left(\frac{s+1}{2}\right)L(s,\chi_d).
\]
\begin{df}\label{kappadef}
For $\varphi\in S(U(\mathbb A_f))$, let 
\[
E(\tau,s;\varphi,+1)=\sum_mA_{\varphi}(s,m,v)\mathbf q^m,
\]
where the Fourier coefficients have Laurent expansions
\[
A_{\varphi}(s,m,v)=b_{\varphi}(m,v)s+O(s^2)
\]
at $s=0$.  For any $\varphi\in S(U(\mathbb A_f))$, define
\[
\kappa_\varphi(m):=\begin{cases} 
\lim_{v\to\infty}b_\varphi(m,v) &\text{if $m>0$,}\\
k_0(0)\varphi(0) &\text{if $m=0$},\\
0 &\text{if $m<0$},
\end{cases}
\]
where 
\begin{equation}\label{kzero}
k_0(0)=\log(d)+2\frac{\Lambda'(1,\chi_d)}{\Lambda(1,\chi_d)}.
\end{equation}
\end{df}
\noindent
For $\varphi=\varphi_\mu=\text{char}(\mu+L)$ we write
\[
A_\mu(s,m,v)=A_{\varphi_\mu}(s,m,v),\; b_\mu(m,v)=b_{\varphi_\mu}(m,v),\;\kappa_\mu(m)=\kappa_{\varphi_\mu}(m).
\]
\begin{theorem}[The $(0,2)$-Theorem]\label{zerotwothm}
Let $F:\h\to S_L\subset S(U(\mathbb A_f))$ be a
weakly holomorphic modular form for $SL_2(\mathbb Z)$ of weight 1, with
Fourier expansion
\[
F(\tau)=\sum_\mu 
F_\mu(\tau)\varphi_\mu=\sum_\mu\sum_mc_\mu(m)\mathbf q^m\varphi_\mu,
\]
where $\mu$ runs over $L^\vee/L$ for some lattice $L$.
Also, assume $c_\mu(m)\in\mathbb Z$ for $m\leq 0$.  Let 
\[
\Phi(z_0,h;F)=\int_{\Gamma\backslash\h}^{\bullet}((\;F(\tau),\theta(\tau,z_0,h)\;))d\mu(\tau).
\]
Then
\[
\int_{SO(U)(\qq)\backslash SO(U)(\mathbb
A_f)}\Phi(z_0,h;F)dh=2\sum_\mu\sum_{m\geq 0} c_\mu(-m)\kappa_\mu(m).
\]
\end{theorem}
\begin{proof}
Our proof is similar to that in \cite{IoBf}.  The integral we want to
compute is given by (\ref{Eonphi}).  Letting $l=-1$ in (\ref{maassonh}), we have
\[
E(\tau,s;\varphi_\mu,-1)v^{-2}=\frac{-4i}{s}\frac{\partial}{\partial\bar{\tau}}\left
\{E(\tau,s;\varphi_\mu,+1)\right\}.
\]
This means we can write
\[
I(s,t)=\frac{1}{2i}\int_{\mathcal F_t}d\Big(\sum_\mu 
F_\mu(\tau)\frac{-4i}{s}E(\tau,s;\varphi_\mu,+1)d\tau\Big).
\]
By Stokes' Theorem, this is
\begin{align}\notag
&=\frac{-2}{s}\int_{\partial\mathcal F_t}\sum_\mu
F_\mu(\tau)E(\tau,s;\varphi_\mu,+1)d\tau\\
&=\frac{-2}{s}\int\limits^{-\frac{1}{2}+it}_{\frac{1}{2}+it}\sum_\mu 
F_\mu(\tau)E(\tau,s;\varphi_\mu,+1)du\notag\\
&=\frac{2}{s}\cdot\textnormal{const. term of}\;\Big(\sum_\mu 
F_\mu(\tau)E(\tau,s;\varphi_\mu,+1)\Big)\bigg|_{v=t}.\label{constterm}
\end{align}
The definition of the regularized integral implies
\begin{align}\notag
&\int_{\Gamma\backslash\h}^{\bullet}((\;
F(\tau),E(\tau,0)\;))d\mu(\tau)=\\
&\underset{\sigma=0}{\textnormal{CT}}\bigg\{\lim_{t\to\infty}\int_{\mathcal 
F_t}\sum_\mu
F_\mu(\tau)E(\tau,0;\varphi_\mu,-1)v^{-\sigma-2}dudv\bigg\}.\notag
\end{align}
\noindent
We need Proposition 2.5 of \cite{IoBf} to hold for $n=0$.  If we use
Proposition 2.6 of \cite{IoBf} and the fact that a factor of $2$ appears in the
Siegel-Weil formula here, then in our notation the analogue of Proposition 2.5
of \cite{IoBf} is
\begin{prop}\label{ssk2.5}
\begin{align}\notag
&\underset{\sigma=0}{\textnormal{CT}}\left\{\lim_{t\to\infty}\int_{\mathcal
F_t}\sum_\mu
F_\mu(\tau)E(\tau,0;\varphi_\mu,-1)v^{-\sigma-2}dudv\right\}\\
&=\lim_{t\to\infty}\left[\int_{\mathcal F_t}\sum_\mu
F_\mu(\tau)E(\tau,0;\varphi_\mu,-1)v^{-2}dudv-2c_0(0)\log(t)\right].\notag
\end{align}
\end{prop}
\begin{proof}
From Lemma \ref{integraltosum}, the left hand side of the desired
identity is
\[
\text{vol}(K)\sum_h\underset{\sigma=0}{\textnormal{CT}}\left\{\lim_{t\to\infty}\int_
{\mathcal
F_t}((F(\tau), \theta(\tau,z_0,h)\;))v^{-\sigma-2}dudv\right\},
\]
where $\text{vol}(K)=\frac{2}{\#(H(\qq)\backslash H(\Af)/K)}$.  Fixing
$h$, we have
\begin{equation}\label{intminusf1}
\underset{\sigma=0}{\textnormal{CT}}\left\{\lim_{t\to\infty}\int_{\mathcal
F_t-\mathcal F_1}((F(\tau),
\theta(\tau,z_0,h)\;))v^{-\sigma-2}dudv\right\}\;+
\int_{\mathcal F_1}((F(\tau),
\theta(\tau,z_0,h)\;))d\mu(\tau).
\end{equation}
The first term in (\ref{intminusf1}) can be written as
\begin{equation}\label{useC}
\underset{\sigma=0}{\textnormal{CT}}\left\{\lim_{t\to\infty}\int\limits_{1}^{t}C(v,h)v^{-\sigma-1}dv\right\},
\end{equation}
where
\begin{align}\notag
C(v,h)&=v^{-1}\int\limits_{-\frac{1}{2}}^{\frac{1}{2}}((F(\tau),
\theta(\tau,z_0,h)\;))du\\
&=\text{const. term of}\;v^{-1}((F(\tau),
\theta(\tau,z_0,h)\;))\notag\\
&=\sum_{\mu}\sum_{\substack{m\in\mathbb Q\\ m\geq 
0}}c_\mu(m)\sum_{\substack{x\in U(\mathbb
Q) \\ Q(x)+m=0}}\varphi_\mu(h^{-1}x)e^{4\pi vQ(x)}.\notag
\end{align}
Then we write (\ref{useC}) as
\begin{equation}\label{breakupC}
\underset{\sigma=0}{\textnormal{CT}}\left\{\lim_{t\to\infty}\int\limits_{1}^{t}[
C(v,h)-c_0(0)]v^{-\sigma-1}dv+\lim_{t\to\infty}\int\limits_{1}^{t}c_0(0)v^{-\sigma-1}dv\right\}.
\end{equation}
As in \cite{IoBf},
\[
\int\limits_{1}^{\infty}[C(v,h)-c_0(0)]v^{-\sigma-1}dv
\]
is a holomorphic function of $\sigma$.  Note, this fact follows, in
part, from Lemma \ref{wt1bound}.  For the other piece
of (\ref{breakupC}) we have
\[
\int\limits_{1}^{t}c_0(0)v^{-\sigma-1}dv=c_0(0)\frac{1}{\sigma}(1-t^{-\sigma}).
\]
This term makes no contribution when we take the limit as $t\to\infty$
followed by the constant term at $\sigma=0$.  We are left with
\[
\lim_{t\to\infty}\left[\int\limits_{1}^{t}
C(v,h)v^{-1}dv-\int\limits_{1}^{t}c_0(0)v^{-1}dv\right]=\lim_{t\to\infty}\left[\int\limits_{1}^{t}
C(v,h)v^{-1}dv-c_0(0)\log(t)\right].
\]
We have the volume
term in front and we sum over $h\in H(\qq)\backslash
H(\Af)/K$, so
this adds on a factor of $2$.  
\end{proof}
\noindent
We point out that the value $c_0(0)$ appearing in (\ref{FexpansionL})
and in Proposition \ref{ssk2.5} is independent of the choice of $L$.  If
we view $F(\tau)\in S(U(\Af))$ as $F(\tau,x)$ for $x\in U(\mathbb
A_f)$, then $c_0(0)$ is the zeroth Fourier coefficient of $F(\tau,0)$.
Proposition \ref{ssk2.5} tells us that
\begin{align}\notag
&\underset{\sigma=0}{\textnormal{CT}}\left\{\lim_{t\to\infty}\int_{\mathcal 
F_t}\sum_\mu 
F_\mu(\tau)E(\tau,0;\varphi_\mu,-1)v^{-\sigma-2}dudv\right\}\\
&=\lim_{t\to\infty}\bigg[\int_{\mathcal
F_t}\sum_\mu
F_\mu(\tau)E(\tau,0;\varphi_\mu,-1)v^{-2}dudv\;-\;2c_0(0)\log(t)\bigg]\notag\\
&=\lim_{t\to\infty}\Big[I(0,t)-2c_0(0)\log(t)\Big].\notag
\end{align}
We need to compute $I(0,t)$.  We have
\begin{equation}\label{Alaurent}
A_{\mu}(s,m,v)=b_{\mu}(m,v)s+O(s^2),
\end{equation}
where there is no
constant term in $A_{\mu}(s,m,v)$ since
$E(\tau,s;\varphi_\mu,+1)$ vanishes at $s=0$.  Then
(\ref{constterm}) implies
\[
I(s,t)=\frac{2}{s}\sum_\mu\sum_mc_\mu(-m)A_{\mu}(s,m,t),
\]
so using (\ref{Alaurent}) we have 
\begin{equation}\label{Iofzero}
I(0,t)=2\sum_\mu\sum_{m} c_\mu(-m)b_\mu(m,t).
\end{equation}
Now we show that parts (i) and (ii) of Proposition 2.11 of \cite{IoBf} hold
for $n=0$.
\begin{prop}\label{prop2.11}
(i) For $m<0$, $b_\mu(m,t)$ decays exponentially as
$t\to\infty$. \\
(ii)
\[
\lim_{t\to\infty}\left(2\sum_\mu\sum_{m<0}c_\mu(-m)b_\mu(m,t)\right)
=0.
\]
\end{prop}
\begin{proof}
If $\varphi_\mu=\otimes_p\varphi_{\mu,p}\in S(U(\Af))$ and  
\[
E(\tau,s;\varphi_\mu,+1)=\sum_mE_m(\tau,s;\varphi_\mu,+1),
\]
then for $m\neq 0$ we have the product formula 
\[
E_m(\tau,s;\varphi_\mu,+1)=A_\mu(s,m,v)\mathbf
q^m=W_{m,\infty}(\tau,s;+1)\prod_p W_{m,p}(s,\varphi_{\mu,p}),
\]
where $W_{m,\infty}(\tau,s;+1)$ and $W_{m,p}(s,\varphi_{\mu,p})$ are the
local Whittaker factors at $\infty$ and $p$, respectively.
Proposition 2.6 (iii) of \cite{KRY} tells us that for $m<0$,
\[
W_{m,\infty}(\tau,0;+1)=0,
\]
and
\[
W'_{m,\infty}(\tau,0;+1)=\pi i\mathbf 
q^m\int\limits_{1}^{\infty}r^{-1}e^{-4\pi|m|vr}dr.
\]
For the finite primes we have
\[
C(m):=\left(\prod_p W_{m,p}(s,\varphi_{\mu,p})\right)\!\bigg|_{s=0}=O(1).
\]
Then
\begin{align}\notag
b_\mu(m,t)&=C(m)W'_{m,\infty}(\tau,0;+1)\\
&=C(m)\pi i\mathbf q^m\int\limits_{1}^{\infty}r^{-1}e^{-4\pi|m|vr}dr,\notag
\end{align}
and we have
\[
\left|b_\mu(m,t)\right|=O\!\left(v^{-1}|m|^{-1}e^{-4\pi|m|v}\right).
\]
This proves (i).  Part (ii) then follows from Lemma \ref{wt1bound}.
\end{proof}
\noindent
Part (ii) of Proposition \ref{prop2.11} tells us that we may ignore the sum on
$m<0$ in (\ref{Iofzero}).  This means our formula for the integral
is
\begin{align}\notag
&\int_{SO(U)(\qq)\backslash SO(U)(\mathbb
A_f)}\Phi(z_0,h;F)dh\;=\\
&\lim_{t\to\infty}\left[2\sum_\mu\sum_{m\geq 0} 
c_\mu(-m)b_\mu(m,t)-2c_0(0)\log(t)\right].\notag
\end{align}
We can improve this by looking at the $m=0$ part. The analogue of
Proposition 2.11 (iii) of \cite{IoBf} is
\begin{lem}\label{kzerozero}
For $m=0$, 
\[
b_0(0,t)-\log(t)=\log(d)+2\frac{\Lambda'(1,\chi_d)}{\Lambda(1,\chi_d)},
\]
and for $\mu\neq 0, b_\mu(0,t)=0$. 
\end{lem}
\begin{proof}
By Theorem 3.1 of \cite{Yang}, we have
\begin{align}\notag
E_0(\tau,s;\varphi_\mu,+1)&=v^{\frac{s}{2}}\varphi_\mu(0)+W_{0,\infty}(\tau,s;+1)\prod_p
W_{0,p}(s,\varphi_{\mu,p})\\
&=v^{\frac{s}{2}}\varphi_\mu(0)-2\pi
i\frac{2^{-s}\Gamma(s)v^{-\frac{s}{2}}}{\Gamma\left(\frac{s}{2}+1\right)\Gamma\left(\frac{s}{2}\right)}\prod_p
W_{0,p}(s,\varphi_{\mu,p}),\notag
\end{align}
which by the duplication formula is
\[
=v^{\frac{s}{2}}\varphi_\mu(0)-\sqrt{\pi}iv^{-\frac{s}{2}}\frac{\Gamma\left(\frac{s+
1}{2}\right)}{\Gamma\left(\frac{s}{2}+1\right)}\prod_p
W_{0,p}(s,\varphi_{\mu,p}).
\]
Theorem 5.2 of \cite{Yang} implies $W_{0,p}(s,\varphi_{\mu,p})=0$ if $\varphi_{\mu,p}$
is not the characteristic function of the local lattice.  So
$b_\mu(0,t)=0$ for $\mu\neq 0$.  Now let
$\mu=0$.  Propositions 2.1 and 6.3 of \cite{Yang} imply
\[
E_0(\tau,s;\varphi_0,+1)=v^{\frac{s}{2}}-\sqrt{\pi}v^{-\frac{s}{2}}\frac{\Gamma\left(\frac{s+
1}{2}\right)L(s,\chi_d)}{\Gamma\left(\frac{s}{2}+1\right)L(s+1,\chi_d)}\mathcal
C_0,
\]
where
\[
\mathcal C_0=2^{\beta_2}\prod_{\substack{q\mid
d \\ q=\text{odd prime}}}q^{-\frac{1}{2}}
\]
and 
\[
\beta_2=\begin{cases}
0 &\text{if $2$ is unramified,}\\
-1 &\text{if $4\mid d$ and $8\nmid d$,}\\
-\frac{3}{2} &\text{if $8\mid d$.}
\end{cases}
\]
Then $\mathcal C_0=d^{-\frac{1}{2}}$.  We have
\begin{align}\notag
E_0(\tau,s;\varphi_0,+1)&=v^{\frac{s}{2}}-v^{-\frac{s}{2}}\frac{\pi^{-\frac{s+1}{2}}\Gamma\left(\frac{s+
1}{2}\right)L(s,\chi_d)}{\pi^{-\frac{s}{2}-1}\Gamma\left(\frac{s}{2}+1\right)L(s+1,\chi_d)}d^{-\frac{1}{2}}\\
&=v^{\frac{s}{2}}-v^{-\frac{s}{2}}\frac{\Lambda(s,\chi_d)}{\Lambda(s+1,\chi_d)}d^{-\frac{1}{2}}.\notag
\end{align}
The functional equation for $\Lambda(s,\chi_d)$ (cf. \cite{Davenport}) is
\[
\Lambda(s,\chi_d)=d^{\frac{1}{2}-s}\Lambda(1-s,\chi_d).
\]
We
normalize $E_0(\tau,s;\varphi_0,+1)$ by
$d^{\frac{s+1}{2}}\Lambda(s+1,\chi_d)$ giving
\begin{align}\notag
E_0^*(\tau,s;\varphi_0,+1)&=d^{\frac{s+1}{2}}v^{\frac{s}{2}}\Lambda(1+s,\chi_d)-d^{\frac{s+1}{2}}v^{-\frac{s}{2}}d^{-s}\Lambda(1-s,\chi_d)\\
&=d^{\frac{s+1}{2}}v^{\frac{s}{2}}\Lambda(1+s,\chi_d)-d^{\frac{1-s}{2}}v^{-\frac{s}{2}}\Lambda(1-s,\chi_d).\notag
\end{align}
Hence,
\begin{align}\notag
E_0^{*,\prime}(\tau,0;\varphi_0,+1)&=2\frac{\partial}{\partial
s}\left\{d^{\frac{s+1}{2}}v^{\frac{s}{2}}\Lambda(1+s,\chi_d)\right\}\Big|_{s=0}\\
&=d^{\frac{1}{2}}\Lambda(1,\chi_d)\left\{\log(d)+\log(v)+2\frac{\Lambda'(1,\chi_d)}{\Lambda(1,\chi_d)}\right\}\notag\\
&=h_k\left\{\log(d)+\log(v)+2\frac{\Lambda'(1,\chi_d)}{\Lambda(1,\chi_d)}\right\}\notag,
\end{align}
by the residue formula.  Then since
$E^{*,\prime}(\tau,0;\varphi_0,+1)=h_kE'(\tau,0;\varphi_0,+1)$, we
have
\[
b_0(0,t)-\log(t)=\log(d)+2\frac{\Lambda'(1,\chi_d)}{\Lambda(1,\chi_d)}.
\]
\end{proof}
\noindent
Now the $m=0$ part is
\[
2\sum_\mu
c_\mu(0)b_\mu(0,t)-2c_0(0)\log(t)=2\sum_{\mu\neq 0}c_\mu(
0)b_\mu(0,t)+2c_0(0)(b_0(0,t)-\log(t)),
\]
and Lemma \ref{kzerozero} tells us that this expression is 
$2c_0(0)k_0(0)$.  
This finishes the proof of Theorem \ref{zerotwothm}.
\end{proof}

\section{Main theorem in the case of signature $(n,2)$}

\subsection{The Rational Splitting $V=V_+\oplus U$}
Now we consider the general case.  Assume that we have a
decomposition $V=V_+\oplus U$ where $V_+$ has signature $(n,0)$ and $U$
has signature $(0,2)$.  For $x\in V$, write $x=x_1+x_2$, $x_1\in V_+, x_2\in U$.
Let $z_0\in D_0$.
Then $R(x,z_0)=-(x_2,x_2)$ so we see 
\[
\varphi_\infty(x,z_0)=e^{-\pi(x,x)_{z_0}}=e^{-\pi[(x_1,x_1)-(x_2,x_2)]}=e^{-\pi(x_1,x_1)}e^{\pi(x_2,x_2)},
\]
which is equal to $\varphi_{\infty,+}(x_1)\varphi_{\infty,-}(x_2)$
for the Gaussians on $V_+$ and $U$, respectively.  We also have
$\omega(g^\prime_\tau)\varphi_\infty=\omega_+(g^\prime_\tau)\varphi_{\infty,+}\otimes\omega_-(g^\prime_\tau)\varphi_{\infty,-}$
for the corresponding Weil representations.  For this decomposition of
$V$, we can write the theta function on $S(V(\mathbb
A_f))$ as a tensor product of two distributions, one on $S(V_+(\mathbb
A_f))$ and one on
$S(U(\mathbb A_f))$.  To see this, let
$\varphi\in S(V(\mathbb A_f))$.  The theta functions are linear, so it
suffices to look at a factorizable Schwartz function
$\varphi=\varphi_+\otimes\varphi_-$.  This gives
\begin{align}\notag
\theta(\tau,z_0,h;\varphi)=\;&v^{-\frac{l}{2}}\sum_{x\in
V(\mathbb
Q)}\omega(g^\prime_\tau)(\varphi_\infty(\cdot,z_0)\otimes\omega(h)\varphi)(x)\\
=\;&v^{-\frac{l}{2}}\sum_{x_1,x_2}(\omega_+(g^\prime_\tau)\varphi_{\infty,+}(x_1)\varphi_+(h_+^{-1}x_1))(\omega_-(g^\prime_\tau)\varphi_{\infty,-}(x_2)\varphi_-(h_-^{-1}x_2))\notag\\
=\;&v^{-\frac{n}{4}}\left(\sum_{x_1}\omega_+(g^\prime_\tau)\varphi_{\infty,+}(x_1)\varphi_+(h_+^{-1}x_1)\right)\notag\\
&\!\!\times\;v^{\frac{1}{2}}\left(\sum_{x_2}\omega_-(g^\prime_\tau)\varphi_{\infty,-}(x_2)\varphi_-(h_-^{-1}x_2)\right)\notag\\
=\;&\theta_+(\tau,z_0,h_+;\varphi_+)\,\theta_-(\tau,z_0,h_-;\varphi_-).\notag
\end{align}
Hence,
\[
\theta(\tau,z_0,h)=\theta_+(\tau,z_0,h_+)\otimes\theta_-(\tau,z_0,h_-),
\] 
where their respective weights
are $\frac{n}{2}$ and $-1$. Since $z_0$ is fixed, we write 
\[
\theta_\pm(\tau,h_\pm)=\theta_\pm(\tau,z_0,h_\pm).
\]

\subsection{The Contraction Map}\label{contract}

Now we describe the main way in which we use the above
factorization of the theta function.  Let $\varphi\in S(V(\mathbb
A_f))$.  Then we can write $\varphi=\sum_j\varphi_+^j\otimes\varphi_-^j$, where
$\varphi_+^j\in S(V_+(\mathbb A_f)), \varphi_-^j\in S(U(\mathbb
A_f))$ and the sum is finite.  We define the \emph{contraction map}
\[
\langle\cdot,\theta_+(\tau,h_+)\rangle:S(V(\mathbb A_f))\to S(U(\mathbb A_f))
\]
by 
\[
\langle\varphi,\theta_+(\tau,h_+)\rangle:=\sum_j\theta_+(\tau,h_+;\varphi_+^j)\,\varphi_-^j.
\]
It is
clear that
\begin{equation}\label{double}
((\;\varphi,\theta(\tau,z_0,h)\;))=((\;\langle\varphi,\theta_+(\tau,h_+)\rangle,\theta_-(\tau,h_-)\;)).
\end{equation}
The expression on the right hand side is nice because it is the pairing of
a function in $S(U(\mathbb A_f))$ and the theta function for $U$.
This is just as in the $n=0$ case.  The value of the contraction map that
we are interested in is $\langle
F(\tau),\theta_+(\tau,1)\rangle$. 
\begin{prop}\label{3facts}
If $F:\h\to S_L$ is a weakly holomorphic modular form of weight
$1-\frac{n}{2}$ and type $\omega$ for $\Gamma'$ whose non-positive
Fourier coefficients lie in $\zz$, then \\
(i) $\langle
F(\tau),\theta_+(\tau,1)\rangle$ is a weakly holomorphic modular form of weight
$1$ and type $\omega_-$ for $\Gamma'$ (cf. Definition \ref{typeomega}), \\
(ii) $\langle
F(\tau),\theta_+(\tau,1)\rangle\in S_{L_-}$ for $L_-=U\cap L$, \\
(iii) The non-positive Fourier coefficients of $\langle
F(\tau),\theta_+(\tau,1)\rangle$ lie in $\zz$.
\end{prop} 
\begin{proof}
By definition,
\begin{equation}\label{bracketMF}
\langle
F(\gamma'\tau),\theta_+(\gamma'\tau,h_+)\rangle=(c\tau+d)\left\langle\omega(\gamma')(F(\tau)),\omega^\vee_+(\gamma')(\theta_+(\tau,h_+))\right\rangle_{
U}.
\end{equation}
Assume that $F(\tau)=\sum_j\varphi_+^j\otimes\varphi_-^j$.  We have 
\[
\omega^\vee_+(\gamma')(\theta_+(\tau,h_+))=\theta_+(\tau,h_+;\omega_+(\gamma')^{-1}\circ\cdot),
\]
so (\ref{bracketMF}) is
\begin{align}\notag
&=(c\tau+d)\left\langle\sum_j\omega_+(\gamma')(\varphi_+^j)\otimes\omega_-(\gamma')(\varphi_-^j),\theta_+(\tau,h_+;\omega_+(\gamma')^{-1}\circ\cdot)\right\rangle_{\!\!U}\\
&=(c\tau+d)\sum_j\theta_+(\tau,h_+;\omega_+(\gamma')^{-1}\omega_+(\gamma')(\varphi_+^j))\omega_-(\gamma')(\varphi_-^j)\notag\\
&=(c\tau+d)\sum_j\theta_+(\tau,h_+;\varphi_+^j)\omega_-(\gamma')(\varphi_-^j)\notag\\
&=(c\tau+d)\omega_-(\gamma')\left(\langle F(\tau),\theta_+(\tau,h_+)\rangle\right).\notag
\end{align}
This proves (i).  
\par In order to
compute the Fourier expansion of $\langle F(\tau),\theta_+(\tau,h_+)\rangle$, we need the expansion of
$\theta_+(\tau,h_+;\varphi_+)$ for $\varphi_+\in S(V_+(\mathbb
A_f))$.  We take $h_+=1$ since the integral we are interested in is
\[
\int_{SO(U)(\mathbb Q)\backslash SO(U)(\mathbb A_f)}\Phi(z_0,h;F)dh.
\]
The explicit $\mathbf q$-expansion of
$\theta_+(\tau,1;\varphi_+)$ is obtained via the action of the Weil representation on
$S(V_+(\mathbb R))$.  In our
particular case,
\begin{align}\notag
\theta_+(\tau,1;\varphi_+)&=v^{-\frac{n}{4}}\sum_{x_1\in V_+(\mathbb
Q)}\omega_+(g^\prime_\tau)\varphi_{\infty,+}(x_1)\varphi_+(x_1)\\
&=v^{-\frac{n}{4}}\sum_{x_1}\omega_+(g^\prime_\tau)e^{-\pi(x_1,x_1)}\varphi_+(x_1),\notag
\end{align}
which by \cite{Shintani} is
\begin{align}\notag
&=v^{-\frac{n}{4}}\sum_{x_1}v^{\frac{n}{4}}e^{2\pi iuQ(x_1)}e^{-\pi v(x_1,x_1)}\varphi_+(x_1)\\
&=\sum_{x_1}e^{2\pi i\tau
Q(x_1)}\varphi_+(x_1)\notag\\
&=\sum_{m\in\mathbb
Q}\Big(\sum_{\substack{x_1\\Q(x_1)=m}}\varphi_+(x_1)\Big)\mathbf q^m.\label{plusexp}
\end{align}
Define
\[
d_{\varphi_+}(m):=\sum_{\substack{x_1\\Q(x_1)=m}}\varphi_+(x_1).
\]
Let $L_+\subset V_+$ be a lattice.  Note that if $\varphi_+$ is the
characteristic function of a coset
$\lambda_++L_+$, then $d_{\varphi_+}(m)$ is an integer which counts the number
of vectors $x_1\in\lambda_++L_+$ such that $Q(x_1)=m$.  Also, $V_+(\mathbb Q)$ is
positive definite so $m\geq 0$ in (\ref{plusexp}).  
\par Now we compute the Fourier
expansion of $\langle F(\tau),\theta_+(\tau,1)\rangle$.  We know
$F(\tau)\in S_L$ for some lattice $L\subset V$.  If we let
$L_+=V_+\cap L$ and $L_-=U\cap L$, then generally the lattice $L$
does not split, i.e., $L\supsetneqq L_++L_-$.  We have
\[
L_++L_-\subset L\subset L^\vee\subset L_+^\vee+L_-^\vee.
\]
Let
\[
L^\vee=\bigcup_\eta(\eta+L),\;L=\bigcup_\lambda(\lambda+L_++L_-),
\]
where $\eta$ and $\lambda$ range over
$L^\vee/L$ and $L/(L_++L_-)$, respectively.  If we write
$\eta=\eta_++\eta_-$ and $\lambda=\lambda_++\lambda_-$, then
\[
L^\vee=\bigcup_\eta\bigcup_\lambda\left(\eta_++\lambda_++L_+\right)+\left(\eta_-+\lambda_-+L_-\right).
\]
Let $F(\tau)=\sum_\eta F_\eta(\tau)\varphi_{\eta+L}$ for
$\varphi_{\eta+L}=\text{char}(\eta+L)$.  Then
\[
\varphi_{\eta+L}=\sum_\lambda\varphi_{\eta_++\lambda_++L_+}\otimes\varphi_{\eta_-+\lambda_-+L_-},
\]
and we have
\[
F(\tau)=\sum_\eta F_\eta(\tau)\sum_\lambda \left(\varphi_{\eta_++\lambda_++L_+}\otimes\varphi_{\eta_-+\lambda_-+L_-}\right).
\]
By definition of the contraction map, this gives
\begin{equation}\label{bigexp}
\langle F(\tau),\theta_+(\tau,1)\rangle=\sum_\eta\sum_\lambda F_\eta(\tau)\theta_+\left(\tau,1;\varphi_{\eta_++\lambda_++L_+}\right)\varphi_{\eta_-+\lambda_-+L_-}.
\end{equation}
From
(\ref{bigexp}), we see that
\[
\langle F(\tau),\theta_+(\tau,1)\rangle\in S_{L_-},
\]
but we point out that
the cosets $\eta_-+\lambda_-+L_-$ need not be incongruent mod
$L_-$.  Let
$c_\eta(m)=c_{\varphi_{\eta+L}}(m)$ and
$d_{\eta_++\lambda_+}(m)=d_{\varphi_{\eta_++\lambda_++L_+}}(m)$.  Then
the Fourier expansion of $\langle
F(\tau),\theta_+(\tau,1)\rangle$ is
\begin{align}\notag
\langle
F(\tau),\theta_+(\tau,1)\rangle&=\sum_\eta\sum_\lambda\Big(\sum_m
c_\eta(m)\mathbf q^m\Big)\Big(\sum_m d_{\eta_++\lambda_+}(m)\mathbf
q^m\Big)\varphi_{\eta_-+\lambda_-+L_-}\\
&=\sum_\eta\sum_\lambda\sum_m\Big(\sum_{m_1+m_2=m}c_\eta(m_1)d_{\eta_++\lambda_+}(m_2)\Big)\mathbf
q^m\varphi_{\eta_-+\lambda_-+L_-}\notag\\
&=\sum_\eta\sum_\lambda\sum_mC_{\eta,\lambda_+}(m)\mathbf q^m\varphi_{\eta_-+\lambda_-+L_-}\notag,
\end{align} 
where we define 
\[
C_{\eta,\lambda_+}(m):=\sum_{m_1+m_2=m}c_\eta(m_1)d_{\eta_++\lambda_+}(m_2).
\]
The coefficients $d_{\eta_++\lambda_+}(m)\in\zz_{\geq 0}$ for $m\geq
0$ and $d_{\eta_++\lambda_+}(m)=0$ if $m<0$.  So assuming
$c_\eta(m)\in\zz$ for $m\leq 0$ implies $C_{\eta,\lambda_+}(m)\in\zz$
for $m\leq 0$, and this finishes the proof of Proposition \ref{3facts}.   
\end{proof}
\noindent
We have seen that the zeroth Fourier coefficient of the modular form
$F$ is
very important.  For example, it gives the weight of $\Psi(F)^2$.  When doing the general case, we use the contraction
map to go from a modular form $F\in S(V(\Af))$ to $\langle
F,\theta_+\rangle\in S(U(\Af))$.  Hence, we will want to know the
zeroth coefficient of $\langle
F,\theta_+\rangle$.  For any modular form $\tilde{F}\in
S(U(\Af))$, define
\[
c_0(0)(\tilde{F})
\]
to be the zeroth Fourier coefficient of $\tilde{F}$.
\begin{cor}
The Fourier expansion of $\langle F(\tau),\theta_+(\tau,1)\rangle$
is
\[
\langle F(\tau),\theta_+(\tau,1)\rangle=\sum_\eta\sum_\lambda\sum_mC_{\eta,\lambda_+}(m)\mathbf q^m\varphi_{\eta_-+\lambda_-+L_-},
\]
where
\[
C_{\eta,\lambda_+}(m)=\sum_{m_1+m_2=m}c_\eta(m_1)d_{\eta_++\lambda_+}(m_2),
\]
and
\begin{equation}\label{c0}
c_0(0)(\langle F,\theta_+\rangle):=c_0(0)(\langle
F(\tau),\theta_+(\tau,1)\rangle)=\sum_\eta\sum_{\substack{\lambda\\\eta_-+\lambda_-=0}}C_{\eta,\lambda_+}(0).
\end{equation}
\end{cor}

\subsection{The $(n,2)$-Theorem}

Recall that the lattice $L$ may not split.  For $\eta\in L^\vee/L$ and
$\lambda\in L/(L_++L_-)$ we write $\eta=\eta_++\eta_-$ and
$\lambda=\lambda_++\lambda_-$.
\begin{df}
Define
\[
\kappa_\eta(m):=\sum_\lambda\sum_{x\in\eta_++\lambda_++L_+}\kappa_{\eta_-+\lambda_-}(m-Q(x)).
\]
\end{df}
\noindent
Note that Definition \ref{kappadef} implies the sum over
$x\in\eta_++\lambda_++L_+$ is finite.

\begin{theorem}[The $(n,2)$-Theorem]\label{mainthm}
Let $F:\h\to S_L\subset S(V(\mathbb A_f))$ be a
weakly holomorphic modular form for $\Gamma'$ of weight $1-\frac{n}{2}$, with
Fourier expansion
\begin{equation}\label{fexp}
F(\tau)=\sum_\eta F_\eta(\tau)\varphi_{\eta}=\sum_\eta\sum_mc_\eta(m)\mathbf q^m\varphi_{\eta},
\end{equation}
where $\varphi_{\eta}=\emph{char}(\eta+L)$ and $\eta$ runs over
$L^\vee/L$.  Also, assume $c_\eta(m)\in\zz$ for $m\leq 0$.  Define 
\[
\Phi(z,h;F):=\int_{\Gamma\backslash\h}^{\bullet}((\;F(\tau),\theta(\tau,z,h)\;))d\mu(\tau).
\]
For $z_0\in D_0$ we have
\begin{align}\notag
&(i)\;\Phi(z_0,h;F)\;\text{is always finite,}\\
&(ii)\;\int_{SO(U)(\mathbb Q)\backslash SO(U)(\mathbb
A_f)}\Phi(z_0,h;F)dh=2\sum_{\eta}\sum_{m\geq 0}c_{\eta}(-m)\kappa_{\eta}(m).\notag
\end{align}
\end{theorem}
\begin{proof}
The regularized integral is given by
\[
\Phi(z_0,h;F)=\int_{\Gamma\backslash\h}^{\bullet}\phi(\tau)d\mu(\tau),
\]
where the integrand is
\begin{align}\notag
\phi(\tau)&=((\;F(\tau),\theta(\tau,z_0,h)\;))\\
&=((\;\langle F(\tau),\theta_+(\tau,1)\rangle,\theta_-(\tau,h_-)\;)),\notag
\end{align}
as in (\ref{double}).  Hence,
\begin{equation}\label{shiftfinite}
\Phi(z_0,h;F)=\Phi(z_0,h_-;\langle F(\tau),\theta_+(\tau,1)\rangle),
\end{equation}
and Proposition \ref{alwaysfinite} implies (\ref{shiftfinite}) is
always finite.  We remark that the regularization process does not
depend on the integrand $\phi(\tau)$.

\par For (ii), using (\ref{shiftfinite}) the desired integral can be written
\begin{align}\notag
&\int_{SO(U)(\mathbb Q)\backslash SO(U)(\mathbb
A_f)}\Phi(z_0,h;F)dh\\
&=\int_{SO(U)(\mathbb Q)\backslash SO(U)(\mathbb
A_f)}\Phi(z_0,h_-;\langle F(\tau),\theta_+(\tau,1)\rangle)dh_-.\label{usezero}
\end{align}
Proposition \ref{3facts} tells us we may apply the $(0,2)$-Theorem to
(\ref{usezero}).  Doing this we see
\begin{align}\notag
(\ref{usezero})&=2\sum_\eta\sum_\lambda\sum_{m\geq
0}C_{\eta,\lambda_+}(-m)\kappa_{\eta_-+\lambda_-}(m)\\
&=2\sum_\eta\sum_\lambda\sum_{m\geq
0}\Big(\sum_{m_1+m_2=-m}c_\eta(m_1)d_{\eta_++\lambda_+}(m_2)\Big)\kappa_{\eta_-+\lambda_-}(m)\notag\\
&=2\sum_\eta\sum_\lambda\sum_{m\geq
0}\Big(\sum_{m_1\leq
0}c_\eta(m_1)d_{\eta_++\lambda_+}(-m-m_1)\Big)\kappa_{\eta_-+\lambda_-}(m)\notag\\
&=2\sum_\eta\sum_\lambda\sum_{m\geq
0}\Big(\sum_{m_1\geq
0}c_\eta(-m_1)d_{\eta_++\lambda_+}(m_1-m)\Big)\kappa_{\eta_-+\lambda_-}(m).\label{bothgeq}
\end{align}
If $m>m_1$, then $d_{\eta_++\lambda_+}(m_1-m)=0$, so
\[
(\ref{bothgeq})=2\sum_\eta\sum_\lambda\sum_{m_1\geq
0}c_\eta(-m_1)\Big(\sum_{0\leq m\leq
m_1}d_{\eta_++\lambda_+}(m_1-m)\kappa_{\eta_-+\lambda_-}(m)\Big).
\]
Then
\begin{align}\notag
&\sum_{0\leq m\leq
m_1}d_{\eta_++\lambda_+}(m_1-m)\kappa_{\eta_-+\lambda_-}(m)\\
&=\sum_{0\leq m\leq
m_1}\left(\#\{x\in\eta_++\lambda_++L_+\mid
Q(x)=m_1-m\}\right)\kappa_{\eta_-+\lambda_-}(m)\notag\\
&=\sum_{\substack{x\in\eta_++\lambda_++L_+\\0\leq Q(x)\leq
m_1}}\kappa_{\eta_-+\lambda_-}(m_1-Q(x))\notag\\
&=\sum_{x\in\eta_++\lambda_++L_+}\kappa_{\eta_-+\lambda_-}(m_1-Q(x)),\notag
\end{align}
since $Q(x)\geq 0$ and $\kappa_{\eta_-+\lambda_-}(m)=0$ for $m<0$.  So 
\[
(\ref{bothgeq})=2\sum_\eta\sum_{m\geq 0} c_\eta(-m)\kappa_\eta(m). 
\]
\end{proof}
We now state an important corollary of Theorem \ref{mainthm}, which
gives the average value of the logarithm of a Borcherds form over CM
points.  As in section \ref{CMsection}, let
$T=\text{GSpin}(U)$ and let $K\subset H(\Af)$ be a compact open
subgroup such that $F:\h\to S_L^K$.  Write $K_T=K\cap T(\Af)$ and
recall that we consider the set of CM points
\[
Z(U)_K=T(\qq)\backslash\Big(D_0\times T(\Af)/K_T\Big)\hookrightarrow X_K.
\]
\begin{cor}\label{withpsi}
(i) When $z_0$ is not in the divisor of the Borcherds form $\Psi(F)$
(i.e., when (\ref{phiandpsi}) holds), the result of Theorem
\ref{mainthm} can be stated as
\[
\sum_{z\in Z(U)_K}\log||\Psi(z;F)||^2=\frac{-2}{\emph{vol}(K_T)}\Big(\sum_{\eta}\sum_{m\geq 0}c_{\eta}(-m)\kappa_{\eta}(m)\Big).
\]
(ii) If $U\simeq k=\qq(\sqrt{-d})$ where $-d$ is an odd fundamental
discriminant, then we have the factorization
\[
\prod_{z\in
Z(U)_K}||\Psi(z;F)||^2={\bf rat}\cdot\left((4d\pi)^{-1} e^{2\frac{L'(0,\chi_d)}{L(0,\chi_d)}}\right)^{\!h_kc_0(0)(\langle
F,\theta_+\rangle)},
\]
where ${\bf rat}\in\qq$ and $c_0(0)(\langle F,\theta_+\rangle)$ is the zeroth Fourier coefficient of $\langle F,\theta_+\rangle$, as defined in (\ref{c0}).  Note that the degree of $Z(U)_K$ is $2h_k$, where $h_k$ is the class number of $k$.  This factorization can also be written as
\[
\prod_{z\in Z(U)_K}||\Psi(z;F)||^2={\bf rat}\cdot\left[(4d\pi)^{-h_k}\prod_{a=1}^{d-1}\Gamma\left(\frac{a}{d}\right)^{w_k\chi_d(a)}\right]^{c_0(0)(\langle
F,\theta_+\rangle)},
\]
where $w_k$ is the number of roots of unity in $k$.  The transcendental factor appearing in this factorization is related to Shimura's period invariants \cite{ShimuraAnn},\cite{GrossCS},\cite{Yoshida}.\\
(iii) 
\[
\log({\bf rat})=-h_k\sum_\eta\sum_{m>0}c_\eta(-m)\Big(\sum_\lambda\sum_{\substack{x\in\eta_++\lambda_++L_+\\Q(x)<m}}\kappa_{\eta_-+\lambda_-}(m-Q(x))\Big).
\]
\end{cor}
\begin{proof}
(i) follows from (\ref{phiandpsi}).  For (ii) and (iii) we have
$\text{vol}(K_T)=\frac{2}{h_k}$, and we will see from Theorem
\ref{kappaformula} of the next section that
\begin{equation}\label{mbiggerzero}
-h_k\sum_\eta\sum_{m>0}c_\eta(-m)\Big(\sum_\lambda\sum_{\substack{x\in\eta_++\lambda_++L_+\\Q(x)<m}}\kappa_{\eta_-+\lambda_-}(m-Q(x))\Big)
\end{equation}
is the logarithm of a rational number ${\bf rat}$.  From
$\Lambda(s,\chi_d)=\pi^{-\frac{s+1}{2}}\Gamma\left(\frac{s+1}{2}\right)L(s,\chi_d)$,
we see
\[
\frac{\Lambda'(1,\chi_d)}{\Lambda(1,\chi_d)}=-\frac{1}{2}\log(\pi)+\Gamma'(1)+\frac{L'(1,\chi_d)}{L(1,\chi_d)}.
\]
So for the corresponding part of (\ref{mbiggerzero}) that involves
$\kappa_0(0)$, we have
\[
-h_kc_0(0)(\langle F,\theta_+\rangle)\left(\log(d)-\log(\pi)+2\Gamma'(1)+2\frac{L'(1,\chi_d)}{L(1,\chi_d)}\right),
\]
which equals
\[
h_kc_0(0)(\langle F,\theta_+\rangle)\left(2\frac{L'(0,\chi_d)}{L(0,\chi_d)}-\log(4d\pi)\right).
\]
The second identity in (ii) follows from the Chowla-Selberg formula
(cf. Proposition 10.10 of \cite{Faltht}), which says
\[
\frac{L'(0,\chi_d)}{L(0,\chi_d)}=\frac{w_k}{2h_k}\sum_{a=1}^{d-1}\chi_d(a)\log
\Gamma\left(\frac{a}{d}\right).
\]
\end{proof}

As an immediate consequence of Corollary \ref{withpsi} and Theorem
\ref{kappaformula} of the next section, we obtain a Gross-Zagier
phenomenon about which primes can occur in the factorization of the
rational part of
\[
\prod_{z\in Z(U)_K}||\Psi(z;F)||^2.
\]
For $F$ as in (\ref{fexp}),
define
\[
m_{\text{max}}=\text{max}\{m>0\mid c_\eta(-m)\neq 0\;\text{for some}\;\eta\}.
\] 
\begin{theorem}\label{smallprimes}
Let $-d$ be an odd fundamental discriminant and assume $U\simeq
k=\qq(\sqrt{-d})$.  Then the only primes which occur in the
factorization of the rational part of
\[
\prod_{z\in Z(U)_K}||\Psi(z;F)||^2
\]
are \\ \\
(i) $q$ such that $q\mid d$, \\ \\
(ii) $p$ inert in $k$ with $p\leq dm_{\emph{max}}$.
\end{theorem}
\noindent
Note that this fact holds for all Borcherds forms and all CM points.  In addition, we point out that the modular form $F$ is not needed in
order to obtain $m_{\text{max}}$.  It can be recovered from the
divisor of $\Psi(F)^2$ (cf. Theorem 1.3 of \cite{IoBf}).

\section{Explicit computation of $\kappa_\mu(t)$ for $t\in\qq_{>0}$}

In order to compute examples of our main theorem, we need
to derive explicit formulas for
$\kappa_\mu(t)$ for $t\in\qq_{>0}$.  Our previous discussion of the
Clifford algebra of $U$ and Lemma \ref{Uisk} imply that, without loss
of generality, we may assume $U=k$ is an imaginary quadratic
field with quadratic form $Q$ given by a negative multiple of the
norm-form.  In this section we assume that $L=\a\subset\O_k$ is an integral ideal and
that $Q(x)=-\frac{Nx}{N\a}$, so that $L^\vee=\mathcal D^{-1}\a$, where
$\mathcal D$ is the different of $k$.  This
is certainly not the most general possible lattice.  Write
$\kappa_\mu(t)$ as $\kappa(t,\mu,\a)$ for $\mu\in\mathcal D^{-1}\a/\a$.  For simplicity, we assume that $k=\qq(\sqrt {-d})$, where $d>3, d\equiv$ 3
(mod 4) and is square-free, so that the prime 2 is not ramified.  Let $\chi$ be the
character of $\qq^\times_{\mathbb A}$ associated to $k$, which is
defined via the global quadratic Hilbert symbol so that
$\chi(t)=(t,-d)_{\mathbb A}$.  Then for a prime $p\leq\infty$, the local character is
$\chi_p(t)=(t,-d)_p$ where $(\;,\;)_p$ is the local quadratic Hilbert symbol.
\par Throughout this section we let $p$ denote an unramified prime and $q$ denote a ramified
prime.  Let $\mu$ be a coset in $\mathcal D^{-1}\a/\a$.  Write $\mu_q$ for the image of $\mu$ under the map
\[
\mathcal D^{-1}\a/\a\to\mathcal D^{-1}\a_q/\a_q,
\]
where $\a_q=\a\otimes_\zz\zz_q$.
For $t\in\qq_{>0}$, we
introduce the function
\[
\rho(t)=\#\{\a\subseteq\O_k\mid N\a=t\}.
\]
This function factors as
\begin{equation}\label{rhoprod}
\rho(t)=\prod_p\rho_p(t),
\end{equation}
where $\rho_p(t)=\rho(p^{\textnormal{ord}_p(t)})$.  The explicit formula for
$\kappa(t,\mu,\a)$ is given by the following theorem.    
\begin{theorem}\label{kappaformula}
For $\mu\in\mathcal D^{-1}\a/\a$ and $t\in\qq_{>0}$,
\[
\kappa(t,\mu,\a)=-\frac{1}{h_k}\prod_{q\mid
d}\emph{char}(Q(\mu_q)+\zz_q)(t)\;\times
\]
\[
\Bigg[\rho(dt)\sum_{\substack{q\mid d\\ \mu_q=0}}\eta_q(t,\mu)(\textnormal{ord}_q(t)+1)\log(q)+\eta_0(t,\mu)\sum_{p\;\textnormal{inert}}(\textnormal{ord}_p(t)+1)\rho\left(\frac{dt}{p}\right)\log(p)\Bigg],
\]
where
\[
\eta_q(t,\mu)=(1-\chi_q(-t))\prod_{\substack{q'\mid d\\ q'\neq q\\ \mu_{q'}=0}}(1+\chi_{q'}(-t))
\]
and
\[
\eta_0(t,\mu)=\prod_{\substack{q\mid d\\ \mu_q=0}}(1+\chi_q(-t)).
\]
We take $\eta_0(t,\mu)=1$ if $\mu_q\neq 0$ for all $q\mid d$.
For $t=0$,
\[
\kappa(0,0,\a)=\log(d)+2\frac{\Lambda'(1,\chi_d)}{\Lambda(1,\chi_d)}.
\]
\end{theorem}
\noindent
Note that when $\mu_q\neq 0$ for all $q\mid d$ we have
$\eta_q(t,\mu)=0$ for all $q$ and $\eta_0(t,\mu)=1$, and so we get a
much simpler formula in this ``generic'' case.
\begin{proof}
The value for $t=0$ is defined in Definition \ref{kappadef}.  For $t>0$,
$\kappa(t,\mu,\a)$ is given by the second term in the Laurent
expansion of a certain Eisenstein series.  These Eisenstein series have
factorizations in terms of local Whittaker functions, and we use
these factorizations to derive the above formula for
$\kappa(t,\mu,\a)$. 
Let $\varphi_{\mu_q}$ be the characteristic function
of the coset $\mu_q$, $X=p^{-s}$, and $\tau=u+iv\in\h$.
Using \cite{Yang} and \cite{KRY}, we have the following
formulas for the normalized local Whittaker functions.  For $\mu=0,$
Lemma 2.3 of \cite{KRY} tells us we only need to consider $t\in\zz,$
and for $t>0$ we have,
\begin{equation}\label{infty}
W^*_{t,\infty}(\tau,s)=\gamma_\infty v^{\frac{1-s}{2}}e(tu)\frac{2i\pi^{\frac{s}{2}}e^{2\pi
tv}}{\Gamma(\frac{s}{2})}\int_{u>2tv}e^{-2\pi u}u^{\frac{s}{2}}(u-2tv)^{\frac{s}{2}-1}du,
\end{equation}
\begin{equation}\label{unr}
W^*_{t,p}(s,\varphi_0)=\sum_{r=0}^{\textnormal{ord}_p(t)}(\chi_p(p)X)^r,
\end{equation}
\begin{equation}\label{ramzero}
W^*_{t,q}(s,\varphi_0)=\gamma_q q^{-\frac{1}{2}}\begin{cases} 
1+(q,-t)_qX^{\textnormal{ord}_q(t)+1} &\text{if ord$_q(t)$ is even,}\\
1+(q,-dt)_qX^{\textnormal{ord}_q(t)+1} &\text{if ord$_q(t)$ is odd}.
\end{cases}
\end{equation}
Here $\gamma_\infty$ and $\gamma_q$ are local factors which do not affect our
global computations since $\gamma_\infty\prod_q\gamma_q=1$, where the product is
over all ramified primes.  For an unramified prime $p$, the local
lattice $\a_p=\a\otimes_\zz\zz_p$ is unimodular.  Note that here is where we
have lost generality by assuming $L=\a$ is an integral ideal.
Since $\a_p$ is unimodular, we only need to consider the Whittaker functions for nonzero
cosets at ramified primes.
For $\mu_q\neq 0$ we have
\begin{equation}\label{ramnonzero}
W^*_{t,q}(s,\varphi_{\mu_q})=\gamma_q q^{-\frac{1}{2}}\text{char}(Q(\mu_q)+\zz_q)(t).
\end{equation}
Note that in (\ref{ramnonzero}), $W^*_{t,q}(s,\varphi_{\mu_q})$ is either a nonzero constant
or is identically zero.  
Following \cite{KRY}, the normalized Eisenstein series has Fourier coefficients given by 
\begin{equation}\label{FCprod}
E_t^*\!\left(\tau,s,\Phi^{1,\mu}\right)=v^{-\frac{1}{2}}d^{\frac{s+1}{2}}W^*_{t,\infty}(\tau,s)\prod_{q\mid
d}W^*_{t,q}(s,\varphi_\mu)\prod_{p\nmid d}W^*_{t,p}(s,\varphi_0).
\end{equation}
Write $t=q^{\alpha_q}u$ where
$\alpha_q=\textnormal{ord}_q(t)$. 
We now show that~(\ref{ramzero}) can be combined into one
nice formula.  
\begin{lem}\label{oneformula}
$W^*_{t,q}(s,\varphi_0)=\gamma_q q^{-\frac{1}{2}}\left(1+\chi_q(-t)X^{\alpha_q+1}\right).$
\end{lem}
\begin{proof}
For $\alpha_q$ even, we have
\[
(q,-t)_q=(-t,q)_q=(-t,-1)_q(-t,-q)_q=(-t,-1)_q(-t,dq)_q\chi_q(-t),
\]
and 
\[
(-t,-1)_q(-t,dq)_q=(-t,-dq^{-1})_q=\left(\frac{-dq^{-1}}{q}\right)^{\!\alpha_q}=1.
\]
For $\alpha_q$ odd,
\begin{align}\notag
(q,-dt)_q&=(-1)^{\frac{q-1}{2}}(q,d)_q(q,t)_q\\
&=(-1,q)_q(q,d)_q(-t,-q)_q(-1,q)_q(-t,-1)_q\notag\\
&=(q,d)_q(-t,-1)_q(-t,dq)_q\chi_q(-t),\notag
\end{align}
and
\begin{align}\notag
(q,d)_q(-t,-1)_q(-t,dq)_q&=(q,d)_q(-t,-dq^{-1})_q\\
&=(-1)^{\frac{q-1}{2}}\left(\frac{dq^{-1}}{q}\right)^{\!\alpha_q}\left(\frac{-dq^{-1}}{q}\right)^{\!\alpha_q}\notag\\
&=1.\notag
\end{align}
So~(\ref{ramzero}) can be rewritten as
\begin{equation}\label{ramzeronew}
W^*_{t,q}(s,\varphi_0)=\gamma_q q^{-\frac{1}{2}}\left(1+\chi_q(-t)X^{\alpha_q+1}\right).
\end{equation}
\end{proof} 
\par Let us first compute $\kappa(t,\mu,\a)$ for $\mu=0$ and
$t\in\mathbb N$.  To do this, we need the following special values for the local
Whittaker functions, cf. Lemma 2.5 and Propositions 2.6 and 2.7 of \cite{KRY}.
\begin{lem}\label{specialvalues} At $s=0$ we have \\
(i) $W^*_{t,\infty}(\tau,0)=-\gamma_\infty
2v^{\frac{1}{2}}e(t\tau)$. \\
(ii) $W^*_{t,p}(0,\varphi_0)=\rho_p(t)$, and if $\rho_p(t)=0$ then
\[
W^{*,\prime}_{t,p}(0,\varphi_0)=\frac{1}{2}(\textnormal{ord}_p(t)+1)\rho_p\left(\frac{t}{p}\right)\log(p).
\]
(iii) $W^*_{t,q}(0,\varphi_0)=\gamma_q q^{-\frac{1}{2}}(1+\chi_q(-t))$,
and if $\chi_q(-t)=-1$ then 
\[
W^{*,\prime}_{t,q}(0,\varphi_0)=\gamma_q
q^{-\frac{1}{2}}(\textnormal{ord}_q(t)+1)\rho_q(t)\log(q).
\]  
\end{lem}
\noindent
Given~(\ref{FCprod}), we consider different cases based on
when one and only one local Whittaker function vanishes at $s=0$.
Since $W^*_{t,\infty}(\tau,0)\neq 0$ for $t\in\mathbb N$, there are
two cases. \\ \\
\underline{Case 1}: $W^*_{t,p}(0,\varphi_0)=0$ for $p$ unramified,
$W^*_{t,p^\prime}(0,\varphi_0)\neq 0\;\forall p^\prime\neq p$. \\ \\
$W^*_{t,p}(0,\varphi_0)=0$ implies that $p$ is inert and ord$_p(t)$ is
odd.  Since $W^*_{t,q}(0,\varphi_0)\neq 0$ for $q$ ramified, we have
$\chi_q(-t)=1$ and 
$W^*_{t,q}(0,\varphi_0)=\gamma_q 2q^{-\frac{1}{2}}$.  Computing the derivative of the
Fourier coefficient we get
\[
E^{*,\prime}_t\!\left(\tau,0,\Phi^{1,0}\right)=W^{*,\prime}_{t,p}(0,\varphi_0)\Bigg[v^{-\frac{1}{2}}d^{\frac{1}{2}}W^*_{t,\infty}(\tau,0)\prod_{q\mid
d}W^*_{t,q}(0,\varphi_0)\prod_{\substack{p^\prime\nmid d \\
p^\prime\neq p}}W^*_{t,p'}(0,\varphi_0)\Bigg]
\]
\[
=\log(p)\frac{1}{2}(\textnormal{ord}_p(t)+1)\rho_p\left(\frac{t}{p}\right)\Bigg[-v^{-\frac{1}{2}}d^{\frac{1}{2}}\gamma_\infty
2v^{\frac{1}{2}}e(t\tau)2^{\nu(d)}\prod_{q\mid d}\gamma_q q^{-\frac{1}{2}}\prod_{\substack{p^\prime\nmid d \\
p^\prime\neq p}}\rho_{p^\prime}(t)\Bigg]
\]
\[
=-\log(p)(\textnormal{ord}_p(t)+1)\rho_p\left(\frac{t}{p}\right)e(t\tau)2^{\nu(d)}\prod_{q\mid d}\rho_q\left(\frac{t}{p}\right)\prod_{\substack{p^\prime\nmid d \\
p^\prime\neq p}}\rho_{p^\prime}\left(\frac{t}{p}\right),
\]
since $\rho_q\left(\frac{t}{p}\right)=1$ and
$\rho_{p^\prime}(t)=\rho_{p^\prime}\left(\frac{t}{p}\right)$, and
where $\nu(d)$ is the number of primes dividing $d$.  So we see
\begin{equation}\label{pzero}
E^{*,\prime}_t\!\left(\tau,0,\Phi^{1,0}\right)=-\log(p)(\textnormal{ord}_p(t)+1)2^{\nu(d)}\rho\left(\frac{t}{p}\right)e(t\tau).
\end{equation}
\underline{Case 2}: $W^*_{t,q}(0,\varphi_0)=0$ for $q$ ramified,
$W^*_{t,p}(0,\varphi_0)\neq 0\;\forall p\neq q$. \\ \\
$W^*_{t,q}(0,\varphi_0)=0$ implies $\chi_q(-t)=-1$ while for any
ramified prime $q^\prime\neq q$ we have $\chi_{q'}(-t)=1$ and
$W^*_{t,q^\prime}(0,\varphi_0)=\gamma_{q^\prime}
2(q^\prime)^{-\frac{1}{2}}$.  In this case, we see
\begin{align}\notag
E^{*,\prime}_t\!\left(\tau,0,\Phi^{1,0}\right)=&W^{*,\prime}_{t,q}(0,\varphi_0)\Bigg[v^{-\frac{1}{2}}d^{\frac{1}{2}}W^*_{t,\infty}(\tau,0)\prod_{\substack{q^\prime\mid d \\
q^\prime\neq q}}W^*_{t,q^\prime}(0,\varphi_0)\prod_{p\nmid d}W^*_{t,p}(0,\varphi_0)\Bigg]\\
=&\gamma_q
q^{-\frac{1}{2}}\log(q)(\textnormal{ord}_q(t)+1)\rho_q(t)\;\times\notag\\
&\Bigg[-v^{-\frac{1}{2}}d^{\frac{1}{2}}\gamma_\infty
2v^{\frac{1}{2}}e(t\tau)2^{{\nu(d)}-1}\prod_{\substack{q^\prime\mid d \\
q^\prime\neq q}}\gamma_{q^\prime}(q^\prime)^{-\frac{1}{2}}\prod_{p\nmid d}\rho_p(t)\Bigg]\notag\\
=&-\log(q)(\textnormal{ord}_q(t)+1)2^{\nu(d)}\rho(t)e(t\tau).\label{qzero}
\end{align}
\par Recall that the definition of $\kappa(t,\mu,\a)$ involves the
non-normalized Eisenstein series, and at $s=0$ we have
$E^{*,\prime}(\tau,0,\Phi^{1,\mu})=h_k E^\prime(\tau,0,\Phi^{1,\mu})$.
This fact and the above analysis, particularly~(\ref{pzero})
and~(\ref{qzero}), give
\[
\kappa(t,0,\a)=
\]
\[
-\frac{2^{\nu(d)}}{h_k}\Bigg(\sum_{q\mid
d}\xi_q(t)(\textnormal{ord}_q(t)+1)\rho(t)\log(q)+\sum_{p\;\textnormal{inert}}\xi_0(t)(\textnormal{ord}_p(t)+1)\rho\left(\frac{t}{p}\right)\log(p)\Bigg),
\] 
where
\[
\xi_q(t)=\begin{cases} 
0 &\text{if $\chi_q(-t)=1$ or $\chi_q(-t)=-1=\chi_{q^\prime}(-t)$, for
some ramified prime}\\
 &\text{ $q^\prime\neq q$,}\\
1 &\text{if $\chi_q(-t)=-1$ and $\chi_{q^\prime}(-t)=1$ for all
ramified primes $q^\prime\neq q$},
\end{cases}
\]
and
\[
\xi_0(t)=\begin{cases} 
0 &\text{if $\chi_q(-t)=-1$ for
some ramified prime $q$,}\\
1 &\text{otherwise}.
\end{cases}
\]
\par Now we compute $\kappa(t,\mu,\a)$ for $\mu\neq 0$.  One main thing to
keep in mind is that there is at least one ramified prime $q$ such that
$\mu_q\neq 0$, but the coset can be zero locally at other ramified
primes.  Write $\mu=(\mu_p)$, where if $p$ is unramified then
$\mu_p=0$ and let $\alpha(\mu)=\#\{q\;\text{ramified}\mid\mu_q=0\}$.
Again, we consider two cases. \\ \\
\underline{Case 1}: $W^*_{t,p}(0,\varphi_0)=0$ for $p$ unramified,
$W^*_{t,p^\prime}(0,\varphi_{\mu_{p^\prime}})\neq 0\;\forall
p^\prime\neq p$. \\ \\
The formula for the derivative of the Fourier coefficient is
\[
E^{*,\prime}_t\!\left(\tau,0,\Phi^{1,\mu}\right)=W^{*,\prime}_{t,p}(0,\varphi_0)\Bigg[v^{-\frac{1}{2}}d^{\frac{1}{2}}W^*_{t,\infty}(\tau,0)\prod_{q\mid
d}W^*_{t,q}(0,\varphi_{\mu_q})\prod_{\substack{p^\prime\nmid d \\
p^\prime\neq p}}W^*_{t,p'}(0,\varphi_0)\Bigg].
\]
Then after cancelling some terms and using Lemma \ref{specialvalues}
and~(\ref{ramnonzero}), we get
\[
=\log(p)\frac{1}{2}(\textnormal{ord}_p(t)+1)\rho_p\left(\frac{t}{p}\right)\Bigg[-
2e(t\tau)2^{\alpha(\mu)}\prod_{\substack{q\mid d \\
\mu_q\neq 0}}\text{char}(Q(\mu_q)+\zz_q)(t)\prod_{\substack{p^\prime\nmid d \\
p^\prime\neq p}}\rho_{p^\prime}(t)\Bigg].
\]
If $q$ is a ramified prime with $\mu_q\neq 0$, then
$W_{t,q}^*(0,\varphi_{\mu_q})\neq 0$ implies ord$_q(t)=-1$.
This means $\rho_q(qt)=1$ and this also equals $\rho_q(dt)$.  If
$\mu_q=0$, then $\rho_q(t)=1=\rho_q(dt)$.  Similarly,
$\rho_p\left(\frac{t}{p}\right)=\rho_p\left(\frac{dt}{p}\right)$ and
$\rho_{p^\prime}(t)=\rho_{p^\prime}(dt)=\rho_{p^\prime}\left(\frac{dt}{p}\right)$.
We also see that if $\mu_q=0$, then
$\text{char}(Q(\mu_q)+\zz_q)(t)=\text{char}(\zz_q)(t)=1$.  So the above
formula is
\begin{equation}\label{pzeromu}
=-2^{\alpha(\mu)}\log(p)(\textnormal{ord}_p(t)+1)\rho\left(\frac{dt}{p}\right)e(t\tau)\prod_{q\mid d}\text{char}(Q(\mu_q)+\zz_q)(t).
\end{equation}
\underline{Case 2}: $W^*_{t,q}(0,\varphi_0)=0$ for $q$ ramified,
$W^*_{t,p}(0,\varphi_{\mu_p})\neq 0\;\forall p\neq q$. \\ \\
Here the derivative is given by
\[
E^{*,\prime}_t\!\left(\tau,0,\Phi^{1,\mu}\right)=W^{*,\prime}_{t,q}(0,\varphi_0)\Bigg[v^{-\frac{1}{2}}d^{\frac{1}{2}}W^*_{t,\infty}(\tau,0)\prod_{\substack{q^\prime\mid d \\
q^\prime\neq q}}W^*_{t,q^\prime}(0,\varphi_{\mu_{q^\prime}})\prod_{p\nmid d}W^*_{t,p}(0,\varphi_0)\Bigg]
\]
\[
=\log(q)(\textnormal{ord}_q(t)+1)\rho_q(t)\Bigg[-
2e(t\tau)2^{\alpha(\mu)-1}\prod_{\substack{q\mid d \\
\mu_q\neq 0}}\text{char}(Q(\mu_q)+\zz_q)(t)\prod_{p\nmid d}\rho_p(t)\Bigg]
\]
\begin{equation}\label{qzeromu}
=-2^{\alpha(\mu)}\log(q)(\textnormal{ord}_q(t)+1)\rho(dt)e(t\tau)\prod_{q\mid d}\text{char}(Q(\mu_q)+\zz_q)(t).
\end{equation}
Note that we do not consider the case where
$W^*_{t,q}(0,\varphi_{\mu_q})=0$ for $\mu_q\neq 0$, since then the
Whittaker function is identically zero and there is no contribution to
the derivative.  Formulas~(\ref{pzeromu}) and~(\ref{qzeromu}) imply that for $\mu\neq 0$,
\[
\kappa(t,\mu,\a)=-\frac{2^{\alpha(\mu)}}{h_k}\prod_{q\mid
d}\text{char}(Q(\mu_q)+\zz_q)(t)\;\times
\]
\begin{equation}\label{formulanonzero}
\Bigg(\sum_{q\mid d}\xi_q(t,\mu)(\textnormal{ord}_q(t)+1)\rho(dt)\log(q)+\sum_{p\;\textnormal{inert}}\xi_0(t,\mu)(\textnormal{ord}_p(t)+1)\rho\left(\frac{dt}{p}\right)\log(p)\Bigg),
\end{equation}
where
\[
\xi_q(t,\mu)=\begin{cases} 
0 &\text{if $\mu_q\neq 0$, or $\mu_q=0$ and $\chi_q(-t)=1$, or
$\chi_q(-t)=-1=\chi_{q^\prime}(-t)$}\\ 
 &\text{for some ramified prime $q^\prime\neq q$ with $\mu_{q^\prime}=0$,}\\
1 &\text{if $\mu_q=0, \chi_q(-t)=-1$, and $\chi_{q^\prime}(-t)=1$ for all
ramified}\\
 &\text{primes $q^\prime\neq q$ with $\mu_{q^\prime}=0$},
\end{cases}
\]
and
\[
\xi_0(t,\mu)=\begin{cases} 
0 &\text{if $\chi_q(-t)=-1$ and $\mu_q=0$ for
some ramified prime $q$,}\\
1 &\text{otherwise}.
\end{cases}
\]
If we take $\mu=0$ in the above equations, we see that
$\xi_q(t,0)=\xi_q(t)$, 
$\xi_0(t,0)=\xi_0(t)$ and $\nu(d)=\alpha(0)$.  Also, when $\mu=0$ then $t\in\mathbb N$ so
$\rho(dt)=\rho(t), \rho\left(\frac{dt}{p}\right)=\rho\left(\frac{t}{p}\right)$ and the
characteristic functions can be ignored.  This means~(\ref{formulanonzero})
holds when $\mu=0$ as well.  We then note that once we sum over $q\mid
d$ with $\mu_q=0$ we can replace $2^{\alpha(\mu)}\xi_q(t,\mu)$ with $\eta_q(t,\mu)$
and we have
\[
\eta_0(t,\mu)=2^{\alpha(\mu)}\xi_0(t,\mu).
\]
This finishes the proof of Theorem \ref{kappaformula}.  
\end{proof}

\bibliographystyle{amsplain}
\bibliography{j3}

\end{document}